\documentclass[12pt]{amsart}

\input{xy} 
\xyoption{all}
\usepackage{graphicx}
\usepackage{color}

\newtheorem{theorem}{Theorem} [section]
\newtheorem{prop}[theorem]{Proposition} 
\newtheorem{lemma}[theorem]{Lemma}
\newtheorem{cor}[theorem]{Corollary}

\numberwithin{equation}{section} 
\numberwithin{figure}{section}

\newcommand\C{{\mathbb C}}
\newcommand\Chat { {\hat{\C}} }

\newcommand\R{{\mathbb R}}
\newcommand\Z{{\mathbb Z}}

\newcommand\N{{\mathbb N}}
\newcommand\D{{\mathbb D}}

\newcommand\cM{\mathcal{M}}
\newcommand\cC{\mathcal{C}}

\newcommand{\cG}{\mathcal {G}}

\newcommand{\cL}{\mathcal{L}}
\newcommand{\cP}{\mathcal{P}}
\newcommand{\cS}{\mathcal{S}}
\newcommand{\cB}{\mathcal{B}}
\newcommand\eps{\varepsilon}

\renewcommand\phi{\varphi}
\newcommand\Aut{\operatorname{Aut}}

\newcommand\iso{\simeq} 

\renewcommand\mod{\operatorname{mod}}  

\newcommand\del{\partial} 

\newcommand\<{\langle} 
\renewcommand\>{\rangle} 

\renewcommand\Im {\operatorname{Im}} 
 
\newcommand\MP {\mathcal{M}}

\newcommand\MM {\mathrm{MM}}
\newcommand\id{\mbox{\rm id}}

\newcommand{\CCC}{\mathcal{C}}

\newcommand{\cl}{\overline}
\newcommand\glue{\mbox{\sc glue}}

\begin{document}

\title{Polynomial basins of infinity}

\author{Laura DeMarco and Kevin M. Pilgrim}

\begin{abstract}  
We study the projection $\pi: \MP_d \to \cB_d$ which sends an affine conjugacy
class of polynomial $f: \C\to\C$ to the holomorphic conjugacy class of the restriction of $f$
to its basin of infinity.  
When $\cB_d$ is equipped with a dynamically natural
Gromov-Hausdorff topology, the map $\pi$ becomes continuous and a homeomorphism on the
shift locus.  
Our main result is that all fibers of $\pi$ are connected. Consequently,
quasiconformal and topological basin-of-infinity conjugacy classes are also connected.  The key ingredient in the proof is an analysis of model surfaces and model maps, branched covers between translation surfaces which model the local behavior of a polynomial.  
\end{abstract}

\date{\today}
\thanks{Research of the first author supported by the National Science Foundation.}

\maketitle

\thispagestyle{empty}

\section{Introduction}

Let $f:\C\to\C$ be a complex polynomial of degree $d\geq 2$.  Iterating $f$ yields a dynamical system.  The plane then decomposes into the disjoint union of  its open, connected {\em basin of infinity} defined by 
	$$X(f) = \{z\in\C: f^n(z) \to \infty \mbox{ as } n\to\infty\}$$
and its complement, the compact {\em filled Julia set} $K(f)$.  

Many naturally defined loci in parameter space (such as the connectedness locus, the shift locus, external rays, their impressions, and parapuzzles) are defined by constraints on the dynamics of $f$ on $X(f)$.  Motivated by this, we study the forgetful map  sending a polynomial $f: \C \to \C$ to its restriction $f: X(f) \to X(f)$ on its basin of infinity.  The basin $X(f)$ is equipped with a dynamically natural translation surface structure.  In this work and its sequels \cite{DP:heights, DP:combinatorics} we exploit this Euclidean perspective to analyze the global structure of moduli spaces of complex polynomials.

\subsection{Connected fibers}
The moduli space $\MP_d$ of complex affine conjugacy classes of degree $d$ polynomials inherits a natural topology from the coefficients of representatives $f$.   Let $\cB_d$ denote the set of conformal conjugacy classes of maps $f: X(f) \to X(f)$, and let $$\pi: \MP_d \to \cB_d$$ be the map sending a polynomial $f$ to its restriction $f|X(f)$.  For each $f \in \MP_d$, the basin of infinity $X(f)$ is equipped with a canonical harmonic Green's function $G_f$ and hence a flat conformal metric $|\del G_f|$ with isolated singularities.  We endow the space $\cB_d$ with the Gromov-Hausdorff topology on the metric spaces $(X(f), 2\,|\del G_f|)$ equipped with the self-map $f: X(f)\to X(f)$; see \S\ref{sec:basins}.  With respect to this topology, the space $\cB_d$ becomes a locally compact Hausdorff metrizable topological space.  The shift locus $\cS_d$ consists of polynomials $f$ for which all $d-1$ critical points lie in $X(f)$.  It forms an open subset of $\MP_d$; its image under $\pi$ is dense in $\cB_d$ (Proposition \ref{prop:shift locus is dense}).   

Recall that a continuous map between topological spaces is {\em monotone} if it has connected fibers. 
Our main result is  

\begin{theorem} \label{pi}
The projection  
$$\pi: \MP_d \to \cB_d$$
is continuous, proper, and monotone.   Furthermore, $\pi$ is a homeomorphism on the shift locus. 
\end{theorem}

\begin{figure}
\includegraphics[width=5in]{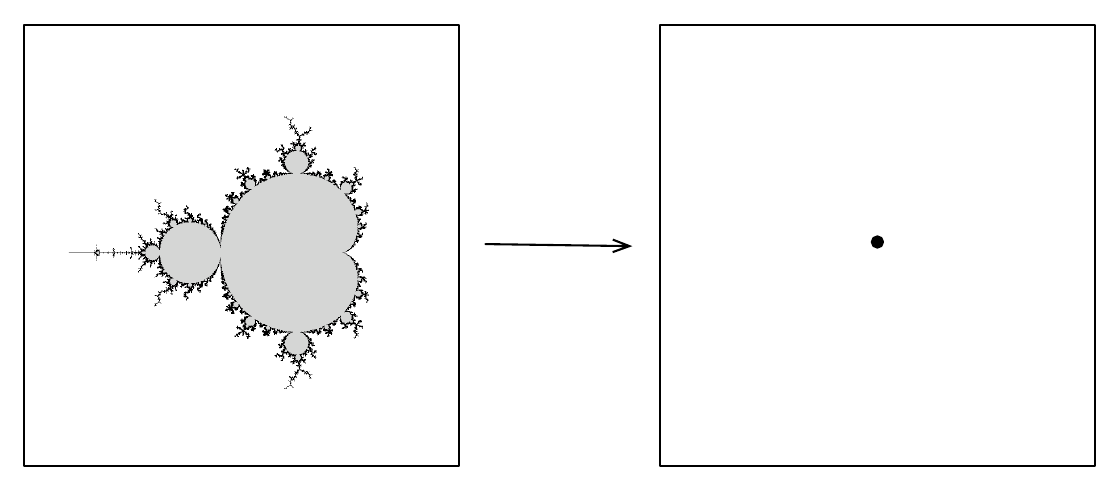}
\caption{For $d=2$, the moduli space $\MP_2$ is isomorphic to the complex plane, as each conjugacy class is uniquely respresented by polynomial $f(z) = z^2 + c$ for $c\in\C$.  The connectedness locus $\CCC_2$ is the much-studied Mandelbrot set.  The projection $\pi: \MP_2 \to \cB_2$ collapses the Mandelbrot set to a point and is one-to-one elsewhere. }
\label{degree 2 projection}
\end{figure}

The key part of Theorem \ref{pi} is the connectedness of fibers, which is already well known in certain important cases.  The fiber over $(z^d, \C\setminus\overline{\D})$ is precisely the {\em connectedness locus} $\mathcal{C}_d$, the set of maps $f$ with connected filled Julia set.  The set $\cC_d$ is known to be cell-like (see \cite [Thm. 8.1]{Douady:Hubbard} for a proof in degree 2,  \cite[Cor. 11.2]{Branner:Hubbard:1} for degree 3, and \cite[Ch. 9]{Lavaurs:thesis} for general degrees), thus connected.  Our theorem gives an alternate proof of its connectedness.  The other extreme is also well known:  for a polynomial in the shift locus, the basin $X(f)$ is a rigid Riemann surface, so such a polynomial is uniquely determined by its basin dynamics.  We exploit this rigidity in the proof of Theorem \ref{pi}.  In the course of the proof we show that the shift locus is connected (Corollary \ref{shift B connected}), a fact which we could not find explicitly stated elsewhere.

\subsection{Topological conjugacy.}
It was observed in \cite[\S8]{McS:QCIII} that any two polynomials $f, g$ which are topologically conjugate on their basins of infinity are in fact quasiconformally conjugate there.  It follows that there is an (analytic) path of polynomials $g_t, 0 \leq t \leq 1$ such that (i) $g_0=g$, (ii) $g_t$ is quasiconformally conjugate on $X(g_t)$ to $g$ on $X(g)$ for all $0 \leq t \leq 1$, and (iii) $g_1$ is conformally conjugate on $X(g_1)$ to $f$ on $X(f)$, i.e. $\pi(g_1)=\pi(f)$.  Since the fiber of $\pi$ containing $f$ is connected, we obtain the following corollary to Theorem \ref{pi}:  

\begin{cor} \label{topological}
Topological or quasiconformal conjugacy classes of basins $(f, X(f))$ are connected in $\MP_d$.
\end{cor}

\subsection{Model maps and sketch proof of Theorem \ref{pi}.}  \label{proof sketch}
Except for the proof of monotonicity, the arguments in the proof are fairly standard.  We record the data of the holomorphic 1-form $\del G_f$ on a basin of infinity $X(f)$.  We use the associated Euclidean structure on the basin to define a Gromov-Hausdorff topology on $\cB_d$.  Continuity and properness of $\pi$ follow from basic properties of $G_f$.  We use the rigidity of basins to deduce that $\pi$ is a homeomorphism on the shift locus.

To treat the monotonicity, we examine the Euclidean structure on a basin of infinity in pieces we call {\em models}:  branched covers between translation surfaces which model the restriction of $f$ to certain subsets of $X(f)$.  We introduce spaces of models, consisting of all branched covers between abstract Riemann surfaces of a special type, and we study the topology of these spaces; via uniformization they may be viewed as subsets of a space of polynomials.  

The idea of the proof of monotonicity in Theorem \ref{pi} is the following.   For each $f$, the Green's function $G_f: X(f) \to (0,\infty)$ is  harmonic and satisfies $G_f(f(z))=d\cdot G_f(z)$.  For $t>0$, let $X_t(f)=\{ z : G_f(z) > t\}$.  Then $f: X_t(f) \to X_{d\cdot t}(f) \subset X_t(f)$, and we may consider the restriction $f|X_t(f)$ up to conformal conjugacy.  

For each $f \in \MP_d$ and $t>0$, we define 
$$\cB(f,t) = \{g\in \MP_d: (g, X_t(g)) \mbox{ is conformally conjugate to } (f, X_t(f))\}.$$
The fiber of $\pi$ containing $f$ may be expressed as the nested intersection  $\bigcap_{t>0}\cB(f,t)$.  We shall show that $\cB(f,t)$ is connected for all (suitably generic) $t$.  The intersection of $\cB(f,t)$ with the shift locus contains a distinguished subset 
$$\cS(f,t) = \{g \in \cB(f,t):  G_g(c) \geq t \mbox{ for all critical points } c \mbox{ of } g \}.$$ 
We show that the space $\cS(f,t)$ is connected by proving that it is homeomorphic to (the finite quotient of) a product of finitely many connected spaces of models.  We construct paths from points in $\cB(f,t)$ to $\cS(f,t)$ by ``pushing up" the critical values; to do this formally, we define the process of {\em gluing} new models into the basin of infinity.  We deduce the connectedness of $\cB(f,t)$ from that of $\cS(f,t)$.    

In fact, the proof of monotonicity of $\pi$ begins like the known proof of connectedness of the  locus $\cC_d$.  When $f$ has connected Julia set, the set $\cB(f,t)$ coincides with 
	$$\cB(t) = \{g \in \MP_d: G_g(c) \leq t \mbox{ for all critical points } c \mbox{ of } g\}$$
for every $t>0$.  It follows from \cite[Cor. 11.2]{Branner:Hubbard:1} and \cite[Ch. 9]{Lavaurs:thesis} that $\cB(t)$ is topologically a closed ball, and its boundary $S(t) = \{g: \max_c G_g(c) = t\}$ is a topological sphere.  The connectedness locus is the nested intersection $\cC_d = \bigcap_t \cB(t)$, showing that $\cC_d$ is cell-like.  By contrast, for general $f$, the topology of $\cB(f,t)$ depends on $f$ and can change as $t$ decreases.

\subsection{Rigidity and other remarks.}  
Intuitively, one might expect that the fibers of $\pi: \MP_d\to\cB_d$ are identified with products of connectedness loci $\cC_{d_i}$ of degrees $d_i \leq d$, each of which is connected.  That is, the affine conjugacy class of a polynomial $f$ should be determined by the conformal conjugacy class of its restriction $(f, X(f))$ together with a finite amount of ``end-data":  the restriction of $f$ to non-trivial periodic components of the filled Julia set $K(f)$.  It is easily seen to hold in degree 2, and it follows in degree 3 by the results of Branner and Hubbard in \cite{Branner:Hubbard:1}, \cite[\S 9]{Branner:Hubbard:2}, where every fiber of $\pi$ in $\MP_3$ is either a point, a copy of the Mandelbrot set $\cC_2$, or the full connectedness locus $\cC_3$.  However, discontinuity of straightening should imply that this intuitive expectation fails in higher degrees; see \cite{Inou:discontinuity}.  

As observed above, for maps $f$ in the shift locus the basin of infinity is a rigid Riemann surface: up to postcomposition with affine maps, there is a unique conformal embedding $X(f) \hookrightarrow \C$.  Consequently, the restricted conformal dynamical system $f: X(f) \to X(f)$ determines the affine conjugacy class of $f: \C \to \C$.   So the projection $\pi: \MP_d\to \cB_d$ is a bijection on the shift locus.  In order for  $X(f)$ to be rigid, the filled Julia set $K(f)$ must be a Cantor set.  Following partial results in \cite[\S 5.4]{Branner:Hubbard:2} and \cite{Emerson:trees}, a converse has been established in \cite{zhai:rigid} and \cite{Yin:Zhai}.  That is, if $K(f)$ is a Cantor set, then the restriction $f: X(f)\to X(f)$ uniquely determines the conformal conjugacy class of $f$.  It follows that the projection $\pi: \MP_d\to \cB_d$ is a bijection on the full Cantor locus in $\MP_d$.  It is not known whether $X(f)$ is always a rigid Riemann surface when $K(f)$ is a Cantor set.

The set $\cS(f,t)$ (introduced in \S\ref{proof sketch}) has been studied by other authors in the special case when all critical points of $f$ have escape rate equal to $t$.  In this case, the set $\cS(f,t)$ is a set $\cG(t)$ independent of $f$; it is the collection of polynomials in $\MP_d$ where all critical points escape at the same rate $t$.  Further, the set $\cG(t)$ is homeomorphic to a finite quotient of the compact, connected space of degree $d$ {\em critical orbit portraits} \cite[Lemma 3.25]{Kiwi:combinatorial}.  The set $\cG(t)$ is equipped with a natural measure $\mu$ inherited from the external angles of critical points.  The Branner-Hubbard {\em stretching} operation deforms a polynomial in the escape locus along a path accumulating on the connectedness locus.  For $\mu$-almost every point of $\cG(t)$, this path has a limit, and the measure $\mu$ pushes forward to the natural bifurcation measure supported in the boundary of the connectedness locus \cite[Thm. 7.20]{Dujardin:Favre:critical}.  

As described above, we use spaces of models to describe the structure of our sets $\cS(f,t)$.  In a sequel to this article, we give alternative descriptions of spaces of models in terms of {\em branched coverings of laminations}.   In this way, spaces of models may be viewed as a generalization of the space of critical portraits.  Using this extra combinatorial structure, we address in \cite{DP:combinatorics} the classification, presently unknown, of the countable set of globally structurally stable conjugacy classes in the shift locus.

\subsection{ Outline.} In section 2 we summarize background from polynomial dynamics.  In section 3 we define the Gromov-Hausdorff topology on $\cB_d$ and prove that the projection $\pi$ is continuous and a homeomorphism on the shift locus.  In section 4 we develop  the theory of model surfaces and model maps.  In section 5, the connectedness of $\cS(f,t)$ is proved, and in section 6 it is applied to complete the proof of Theorem \ref{pi}.

\subsection{ Acknowledgement.}  We would like to thank Curt McMullen for his helpful suggestions, and Hiroyuki Inou, Chris Judge, and Yin Yongcheng  for useful conversations.  We are especially grateful to the anonymous referee for numerous and detailed comments.


\bigskip\bigskip\section{Spaces of polynomials}

In this section we introduce the moduli spaces $\MP_d$ and give some background on polynomial dynamics.

\subsection{Polynomial dynamics.}  
Let $f$ be a complex polynomial of degree $d \geq 2$.  The filled Julia set  
\[ K(f) = \{ z : \mbox{ the sequence } f^n(z) \mbox{ is bounded} \}\]
is compact, and its complement $X(f) = \C \setminus K(f)$ is open and connected.  For $t \in [0,\infty)$ define $\log^+(t)=\max\{0, \log t\}$.  The function 
\[ G_f: \C \to [0,\infty)\]
given by 
\[ G_f(z) = \lim_{n \to \infty} \frac{1}{d^n}\log^+|f^n(z)| \]
measures the rate at which the point $z$ escapes to infinity under iteration of $f$.   It vanishes exactly on $K(f)$, is harmonic on $X(f)$, and on all of $\C$ it is continuous, subharmonic, and satisfies the functional equation $G_f(f(z)) = d \cdot G_f(z)$ (see e.g.~\cite[\S 18]{milnor:dynamics}).  

For $t>0$, we define 
\[ X_t(f) = \{z : G_f(z)>t\}, \;\;\; Y_t(f) = \{ z : t \leq G_f(z) \leq 1/t\},\]
and we set 
\[ M(f) = \max \{ G_f(c) : f'(c)=0\}, \;\;\; m(f) = \min\{ G_f(c) : f'(c)=0\}.\]

\subsection{Monic and centered polynomials}
Every polynomial 
	$$f(z) = a_dz^d + a_{d-1} z^{d-1} + \cdots + a_0,$$
with $a_i\in \C$ and $a_d\not=0$, is conjugate by an affine transformation $A(z) = az+b$ to a polynomial which is monic ($a_d=1$) and centered ($a_{d-1} = 0$).  The monic and centered representative is not unique, as the space of such polynomials is invariant under conjugation by $A(z) = \zeta z$ where $\zeta^{d-1} = 1$.  In this way, we obtain a finite branched covering
	$$\cP_d \to \MP_d$$
from the space $\cP_d \iso\C^{d-1}$ of monic and centered polynomials to the moduli space $\MP_d$ of conformal conjugacy classes of polynomials.  Thus, $\MP_d$ has the structure of a complex orbifold of dimension $d-1$.  The functions $M(f)$ and $m(f)$ are continuous and invariant under affine conjugation, and $M(f)$ is proper  \cite[Prop. 3.6]{Branner:Hubbard:1}.  As a function of $f$ and of $z$, $G_f(z)$ can be expressed as a locally uniform limit of the pluriharmonic functions $\log|f^n(z)|$ where $G_f(z) >  0$; therefore, $G_f(z)$ is pluriharmonic on the locus in $\cP_d$ where $G_f(z)>0$; see the proof of \cite[Prop. 1.2]{Branner:Hubbard:1}.

It is sometimes convenient to work in a space with marked critical points.  Let $\mathcal{H} \subset \C^{d-1}$ denote the hyperplane given by $\{(c_1, \ldots, c_{d-1}) : c_1 + \ldots + c_{d-1}= 0\}$.  Then the map 
\[ \rho: \mathcal{H} \times \C \to \cP_d\]
given by 
\begin{equation}
\label{eqn:param}
 \rho(c_1, \ldots, c_{d-1}; a)(z) = \int_0^z d\cdot \prod_{i=1}^{d-1}(\zeta - c_i)\; d\zeta + a 
 \end{equation}
gives a polynomial parameterization of $\cP_d$ by the location of the critical points and the image of the origin.  Setting $$\cP^\times_d=\mathcal{H}\times \C,$$ we refer to $\cP^\times_d$ as the space of {\em critically marked} polynomials.

\subsection{External rays and angles.} 
Fix a monic and centered polynomial $f$ of degree $d>1$.  Near infinity, there is a conformal change of coordinates which conjugates $f$ to $z\mapsto z^d$.  The local conjugating isomorphism is unique up to multiplication by a $(d-1)$st root of unity and is therefore uniquely determined if required to have derivative 1 at infinity (see e.g. \cite[Prop. 1.4]{Branner:Hubbard:1}).  It extends to an isomorphism
	$$\phi_f: \{G_f > M(f)\} \to \{z\in\C: |z| > e^{M(f)}\}$$
called the {\em B\"ottcher map}, satisfying $G_f = \log |\phi_f|$.
For each fixed $\theta\in \R/2\pi \Z$, the preimage under  $\phi_f$ of the ray $\{re^{i\theta}: r>e^{M(f)}\}$ is called the {\em external ray of angle $\theta$} for $f$.  There are exactly $d-1$ {\em fixed} external rays mapped to themselves under $f$; their arguments are asymptotic to $2\pi k/(d-1)$ near infinity for $k=0, \ldots, d-2$.  

On $\{|z|>e^{M(f)}\}$, for each angle $\theta$, the external ray of angle $\theta$ coincides with a gradient flow line of $G_f$.   This ray can be extended uniquely to all radii $r>1$ provided that when flowing downward, the trajectory does not meet any of the critical points of $G_f$, i.e. critical points of $f$ or any of their iterated inverse images.  It follows that for all but countably many $\theta$, the external ray of angle $\theta$ admits such an extension, i.e. is {\em nonsingular}.  We see then that the external rays of $f$ define a singular {\em vertical foliation} on $X(f)$ which is orthogonal to the singular {\em horizontal foliation} defined by the level sets of $G_f$.  These foliations coincide with the vertical and horizontal foliations associated to the holomorphic 1-form
	$$\omega_f = 2i \, \del G_f$$
on $X(f)$.  We will exploit this point of view further in the next section.  

We emphasize that, by definition, $\MP_d$ is a quotient of $\cP_d$ by the cyclic group of order $d-1$ acting by conjugation via rotations of the plane centered at the origin.  Therefore, given an element of $\MP_d$, it defines a conjugacy class of a dynamical system on a Riemann surface isomorphic to the plane, and it defines an identification of this surface with the plane, up to this rotational ambiguity.   Thus, given an element of $\MP_d$, together with a choice of fixed external ray, there is a unique such identification sending this chosen fixed external ray to the external ray whose asymptotic argument is zero.

\subsection{Critical values}
Throughout this work, we will make implicit use of the following fact, which is easily proved using the Riemann-Hurwitz formula.  If $f: \C \to \C$ is a nonconstant polynomial and $D \subset \C$ is a bounded Jordan domain whose boundary is disjoint from the critical values of $f$, then $f^{-1}(D)$ is a finite union of Jordan domains, the restriction $f: f^{-1}(D) \to D$ is a proper branched covering of degree $\deg(f)$, and $f^{-1}(D)$ consists of a single component if and only if $D$ contains all critical values of $f$.

The following two lemmas have nothing to do with dynamics and will be used in the proof of the connectedness and compactness of the space of local models (\S\ref{spaces of local models}).   The proof of Lemma \ref{cvmap} is a non-dynamical version of the proof in \cite[Prop. 3.6]{Branner:Hubbard:1} showing properness of $f\mapsto M(f)$.

\begin{lemma} \label{cvmap}
Let $\widetilde{\nu}: \cP^\times_d \to \C^{d-1}$ be the map sending a critically marked polynomial to its ordered list of critical values and $\nu: \cP_d \to \C^{d-1}/S_{d-1} \iso \C^{d-1}$ the map sending a polynomial to its unordered set of critical values.  Then $\widetilde{\nu}$ and $\nu$ are proper.  Moreover, $\widetilde{\nu}$ is a polynomial map, and has the property that any path in the codomain can be lifted (not necessarily uniquely) to a path in the domain.
\end{lemma}

\proof   Equation (\ref{eqn:param}) shows the map $\widetilde{\nu}$ is polynomial.   Fix $R>0$.  Suppose $f(z)=z^d + a_{d-2}z^{d-2}+\ldots + a_1 z + a_0$ belongs to $\cP^\times_d$ and the critical values of $f$ lie in $D_R$.  There is a unique univalent analytic map $\psi_f: \C\setminus D_{R^{1/d}} \to \C$ tangent to the identity at infinity and satisfying $f\circ \psi_f(w)=w^d$.  By the Koebe $1/4$-theorem, $\psi_f(\C\setminus D_{R^{1/d}}) \supset \C \setminus D_{4R^{1/d}}$, so the critical points of $f$ are contained in $D_{4R^{1/d}}$.  It follows that the coefficients $a_{d-2}, \ldots, a_1$ are bounded in modulus by a constant $C_1(R)$.  Since in addition $f$ is assumed monic,  the map $f$ is Lipschitz on $D_{4R^{1/d}}$ with constant $C_2(R)$, so the image $f(D_{4R^{1/d}})$ has diameter less than a constant $C_3(R)$.   Since the critical values of $f$ lie in $D_R$, the image $f(D_{4R^{1/d}})$ meets $D_R$, and so $|a_0|=|f(0)|$ is bounded by a constant $C_4(R)$ as well.  Hence $\widetilde{\nu}$ and $\nu$ are proper.  The final assertion about path lifting is well-known; see \cite[\S III.B]{Gunning:Rossi}.
\qed

\begin{lemma}  \label{compactconnected}
Let $C$ be any compact and path-connected subset of $\C$.  The subset of $\cP_d$ with all critical values in $C$ is compact and path-connected.
\end{lemma}

\proof  
Consider the diagram 
$$\xymatrix{ \cP^\times_d \ar[d]_\rho \ar[r]^{\widetilde{\nu}} & \C^{d-1} \ar[d] \\
		\cP_d   \ar[r]^\nu & \C^{d-1}/ S_{d-1} }$$
The right-hand vertical map is proper, and Lemma \ref{cvmap} implies the horizontal maps are proper, so the compactness conclusion holds.   

Fix now $v\in C$.  
There is a unique monic and centered polynomial with a single critical value at $v$ of multiplicity $d-1$.  It is $f_1(z) = z^d + v$.  For any other $f$ with all critical values in $C$, we can construct a path to $f_1$.  Let $(v_1, \ldots, v_{d-1})\in \C^{d-1}$ be a labelling of the critical values of $f$ (listed with multiplicity).  Choose a continuous deformation of these points $v_i(t) \in C$ for $t\in [0,1]$ so that 
\begin{itemize}
\item[(i)] $v_i(0) = v_i$ for all $i$, 
\item[(ii)] $v_i(1) = v$ for all $i$.  
\end{itemize}
By Lemma \ref{cvmap} the motion of labelled critical values can be  lifted under $\widetilde{\nu}$; it can then projected under $\rho$ to obtain a path from $f$ to $f_1$ for which the corresponding maps all have critical values lying in $C$.  Hence this set of polynomials is path-connected.
\qed


\bigskip\bigskip\section{Restricting to the basin of infinity}  
\label{sec:basins}

Recall that $\cB_d$ denotes the space of conformal conjugacy classes of $(f, X(f))$.  Here we introduce the {\em Gromov-Hausdorff topology} on $\cB_d$ and begin the analysis of the restriction map 
	$$\pi: \MP_d \to \cB_d.$$
We prove the continuity of $\pi$ and show that it is a homeomorphism on the shift locus in $\MP_d$.  

\subsection{The conformal metric $|\omega|$ on the basin of infinity} \label{metric}
Fix a polynomial $f$ of degree $d \geq 2$.  On its basin of infinity $X(f)$, recall that $G_f$ denotes the harmonic escape rate function  and 
	$$\omega_f=2i\,\partial G_f$$ 
the corresponding holomorphic $1$-form, so that $|\omega_f|$ is the associated singular flat conformal metric.  In this way, the pair $(X(f), \omega_f)$ becomes a horizontal translation surface with height function $G_f$.  Note that the height $G_f(z)$ of any point $z\in X(f)$ coincides with its $|\omega_f|$-distance to the lower ends of $X(f)$.  Recall that 
	$$M(f) = \max\{G_f(c): f'(c)=0\}$$ 
denotes the maximal critical height of $f$ and 
	$$m(f) = \min\{G_f(c): f'(c) = 0\}$$
the minimal critical height.  

The zeros of $\omega_f$ coincide with the critical points of $f$ in $X(f)$ and all of their preimages by the iterates $f^n$.  The neighborhood $\{z: G_f(z)> M(f)\}$ of infinity is isometric to a half-infinite Euclidean cylinder of radius 1.  In fact,  if $L$ is any horizontal leaf of $\omega_f$ at height $c$, and if the {\em level} of $L$ is defined as the integer 
\begin{equation} \label{level}
	l(L) = \min\{n\geq 0: d^nc \geq M(f)\},
\end{equation}
then the length of $L$ is
	$$\int_L |\omega_f|  = \frac{\deg (f^{l(L)}|L)}{d^{l(L)}} \; 2\pi.$$
Further, if $A$ is a connected component of $\{a < G_f < b\}$ which is topologically an annulus, then it is isometric to a cylinder of height $b-a$ and circumference $\int_L |\omega_f|$ for any horizontal leaf $L$ in $A$.

On the basin of infinity, conformal and isometric conjugacies are the same thing.

\begin{lemma}  \label{isometric=conformal}
Two polynomials $f$ and $g$ are conformally conjugate on their basins of infinity if and only if they are isometrically conjugate with respect to the conformal metrics $|\omega_f|$ and $|\omega_g|$. In particular, the escape rates of the critical points are isometric conjugacy invariants.
\end{lemma}

\proof
If $f$ and $g$ are conformally conjugate, then the conjugacy sends $\omega_f$ to $\omega_g$, and therefore their basins  are isometrically conjugate.  Conversely, the conformal metric $|\omega_f|$ determines the complex structure on $X(f)$, so an isometry must be a conformal isomorphism.  
\qed

\subsection{The topology of $\cB_d$.} \label{basins:topology}
We define here the {\em Gromov-Hausdorff topology} on the space of triples $(f, X(f), |\omega_f|)$ in $\cB_d$.   We will see that this topology is fine enough to guarantee that $\cB_d$ contains a homeomorphic copy of the shift locus $\cS_d$, but is also coarse enough to guarantee that many natural maps and operations on $\cB_d$ are continuous.

Given a point of $\cB_d$ represented by $(X(f), f, |\omega_f|)$, a neighborhood basis consists of the collection of open sets $U_{t, \epsilon}(f) \subset \cB_d$, where $\epsilon > 0$ and $1/t > M(f)$, 
defined as follows.  Recall that 
	$$Y_t(f) = \{ t \leq G_f(z) \leq 1/t\}.$$  
Let $\rho_f(\cdot,\cdot)$ denote the distance function on $Y_t(f)$ induced by the conformal metric $|\omega_f|$ (since $Y_t(f)$ is convex, $\rho_f$ coincides with the restriction of the length metric to $Y_t(f)$).  Now suppose $g$ is a polynomial with $1/t > M(g)$.  An $\eps$-conjugacy between $f|Y_t(f)$ and $g|Y_t(g)$ is a relation which is nearly the graph of an isometric conjugacy.  That is, as a relation, it is a subset $\Gamma \subset Y_t(f)\times Y_t(g)$ such that 
\begin{enumerate}
\item {\bf nearly surjective:} 
\begin{enumerate}
\item for every $a\in Y_t(f)$, there exists a pair $(x,y)\in \Gamma$ such that $\rho_f(a,x) < \eps$,
\item for every $b\in Y_t(g)$, there exists a pair $(x,y)\in \Gamma$ such that $\rho_g(b,y) < \eps$,
\end{enumerate}
\item {\bf nearly isometric:} if $(x,y)$ and $(x',y')$ are in $\Gamma$, then 
	$$|\rho_f(x,x') - \rho_g(y,y')| < \eps,$$
	and
\item {\bf nearly conjugacy:} for each $(x,y)\in \Gamma$ such that $(f(x), g(y))$ lies in $Y_t(f)\times Y_t(g)$, there exists $(x',y')\in \Gamma$ such that $\rho_f(x', f(x)) <\eps$ and $\rho_g(y', g(y)) < \eps$.
\end{enumerate}
The set $U_{t,\eps}(f)$ consists of all triples $(g, X(g), |\omega_g|)$ for which $1/t > M(g)$ and for which there is an $\eps$-conjugacy between $f|Y_t(f)$ and $g|Y_t(g)$. 
\medskip

\noindent{\bf Remark:}  If $1/t > M(f)$ then by B\"ottcher's theorem the restriction $f|_{\{G_f > 1/t\}}$ is holomorphically conjugate to the restriction of $z \mapsto z^d$  acting on $\{\log |z| > 1/t\}$.  Hence if $f, g$ are any two polynomials of degree $d$ and $1/t>\max\{M(f), M(g)\}$, then the restrictions $f|_{\{G_f>1/t\}}$ and $g|_{\{G_g>1/t\}}$ are holomorphically conjugate.  The set of such conjugacies is naturally identified with the group of isometric automorphisms of $f|_{\{G_f>1/t\}}$,  which is the group generated by the rigid rotation of order $d-1$.  It follows that if in addition $f, g$ are $\epsilon$-conjugate via a relation $\Gamma \subset Y_t(f) \times Y_t(g)$, then there is an extension of $\Gamma$ to a relation on $X_t(f) \times X_t(g)$ which gives an  $\epsilon$-conjugacy from $f$ to $g$: an $\epsilon$-conjugacy must send a point on $\{G_f=1/t\}$ with external angle $\theta$ that is fixed under $t \mapsto d\cdot t \bmod 2\pi$ to a point which is $\epsilon$-close to a point $\theta'$ whose external angle is also so fixed; the extension is given by the unique rotation sending $\theta$ to $\theta'$.  
\medskip 

For $t>0$ let 
\[ \cB_{d,t} = \{ (X_t(f), f, |\omega_f|) : f \in \MP_d \} / \sim \]
where $f \sim g$ if there is a holomorphic isomorphism $h: X_t(f) \to X_t(g)$ such that $h\circ f = g \circ h$; we denote by $\pi_t: \MP_d \to \cB_{d,t}$ the corresponding projection.  

We equip $\cB_{d,t}$ with the analogous Gromov-Hausdorff topology: given $f$, a neighborhood basis is given by sets $V_{s, \eps}$, where  $\epsilon>0$ and $1/s > \max\{1/t, M(f)\}$; the set $V_{s, \eps}$ consists of all triples $(g, X_t(g), |\omega_g|)$ for which there is an $\epsilon$-conjugacy from $f|X_t(f)$ to $g|X_t(g)$.  By construction, the projection $\pi_t$ factors as a composition of $\pi$ with the natural projection $\cB\to\cB_{d,t}$.  

\begin{lemma} \label{continuous}
The projection $\pi$ is continuous, surjective, and proper; the projection $\cB_d \to \cB_{d,t}$ induced by $\pi_t$ is continuous and surjective.  
\end{lemma}

\proof Surjectivity holds by definition.  
If $f_k\to f$ in $\MP_d$, then there are polynomial representatives which converge uniformly on compact subsets of $\C$, and the escape-rate functions $G_{f_k}$ converge to $G_f$ by \cite[Proposition 1.2]{Branner:Hubbard:1}.    Fix $t>0$ with $1/t > M(f)$, so that $1/t > M(f_k)$ for all sufficiently large $k$.  The compact sets $Y_t(f_k)$ converge to $Y_t(f)$ in the Hausdorff topology on compact subsets of $\C$, and the action of $f_k$ on $Y_t(f_k)$ converges to that of $f$ on $Y_t(f)$ (with respect to the Euclidean metric on $\C$).  Moreover, the escape-rate functions $G_{f_k}$ and $G_f$ are harmonic near the sets $Y_t(f), Y_t(f_k)$, so uniform convergence implies also the convergence of their derivatives.   Therefore, the 1-forms $\omega_{f_k}$ on $Y_t(f_k)$ converge to $\omega_f$ on $Y_t(f)$ and so the conformal metrics $|\omega_{f_k}|$ on $Y_t(f_k)$ converge to the conformal metric $|\omega_f|$ on $Y_t(f)$.  More precisely:  let $\Gamma_k$ be the graph of the identity map on $Y_t(f_k)\cap Y_t(f)$, regarded as a relation on $Y_t(f_k)\times Y_t(f)$.  For all large enough $k$, the graph $\Gamma_k$ defines an $\eps$-conjugacy between $f_k|Y_t(f_k)$ and $f|Y_t(f)$, and as remarked above it extends to an $\eps$-conjugacy between $f_k|X_t(f_k)$ and $f|X_t(f)$.  Therefore $\pi$ and $\pi_t$ are continuous.
Properness of the map $\pi$ follows from the known fact  that $f \mapsto M(f)$ is proper \cite[Prop. 3.6]{Branner:Hubbard:1}.   
\qed

\begin{lemma}\label{lemma:topology}
The spaces $\cB_d$ and $\cB_{d,t}$ equipped with the Gromov-Hausdorff topology are Hausdorff, locally compact, second-countable, and metrizable.  Moreover, $\cB_d$ is homeomorphic to the quotient space of $\MP_d$ obtained by identifying the fibers of $\pi$ to points.  
\end{lemma}

\proof A standard application of the definitions and a diagonalization argument shows that $\cB_d$ is Hausdorff.  By definition, the topology is first-countable.  By \cite[Thm. 5, p. 16]{Daverman:decompositions} it follows that the Gromov-Hausdorff and quotient topologies coincide.  Local compactness follows from continuity and properness of the projection $\pi$.   Metrizability follows from \cite[Prop. 2, p. 13]{Daverman:decompositions} and second-countability follows.  
\qed
\medskip

\subsection{The shift locus and rigid Riemann surfaces.}  \label{basins:shift}
Recall that the shift locus $\mathcal{S}_d$ is the collection of polynomials in $\MP_d$ where all critical points escape to $\infty$ under iteration.

A planar Riemann surface $X$ is {\em rigid} if a holomorphic embedding $X\hookrightarrow \Chat$ is unique up to postcomposition with conformal automorphisms of $\Chat$.  Equivalently, the complement of $X$ in $\Chat$ has  absolute area zero; that is, the spherical area of $\Chat\setminus X$ is 0 under any holomorphic embedding.  Further, an absolute area zero subset of the plane is removable for locally bounded holomorphic functions with finite Dirichlet integral \cite[\S IV.4]{Ahlfors:Sario}.  In particular, if the complement of $X\subset \C$ has absolute area zero, then any proper holomorphic degree $d$ self-map $X \to X$ extends uniquely to a degree $d$ rational function $\Chat \to \Chat$.  In \cite[\S2.8]{McMullen:CDR}, McMullen showed that an open subset $X\subset \C$ is rigid if it satisfies the {\em infinite-modulus condition}:   for each $n\in\N$, there is a finite union of disjoint unnested annuli $E_n\subset X$, contained in the bounded components of $\C\setminus E_{n-1}$, such that for each nested sequence of connected components $\{A_n\subset E_n\}_n$, we have $\sum_n \mod A_n = \infty$, and the nested intersection of the bounded components of $\Chat\setminus E_n$ is precisely $\C\setminus X$.

Though it has not previously been stated in quite this way, the proof of the following lemma is well-known (see e.g. \cite[Lemma 3.2]{Blanchard:Devaney:Keen}, \cite[\S 5.4]{Branner:Hubbard:2},  \cite[Remark, p. 423]{Branner:cubics}, \cite[\S 2.8]{McMullen:CDR}).

\begin{lemma}  \label{homeo}
The projection $\pi$ is a homeomorphism on the shift locus $\cS_d$.
\end{lemma}

\proof 
When $f$ is in the shift locus, it is easy to see that the basin $X(f)$ satisfies the infinite-modulus condition.  Consider an annulus $A=\{a < G_f(z) < b\}$  with $M(f) < a < b < d\cdot M(f)$ and disjoint from the critical orbits.  Since $f$ is in the shift locus, the iterated preimages  of this annulus map with uniformly bounded degree onto $A$. Hence each such preimage has modulus at least $m>0$.  It follows that there is a unique embedding of $X(f)$ into $\Chat$, up to an affine transformation, sending infinity to infinity.  Furthermore, the complement $\Chat\setminus X(f)$ is removable for $f$, so $f: X(f) \to X(f)$ extends uniquely to a rational map $f: \Chat\to\Chat$ which is totally ramified at infinity.  In other words, up to affine conjugacy, we can reconstruct the polynomial $f: \C \to \C$ from its restriction $f: X(f) \to X(f)$.  

Lemma \ref{continuous} implies that $\pi$ is a continuous bijection on the shift locus.  The image is Hausdorff, and the domain is locally compact.  It follows that the map $\pi$ is a local homeomorphism and therefore a global homeomorphism since $\pi$ is proper. 
\qed


\bigskip\bigskip\section{Models of a polynomial branched cover}  \label{sec:local models}

In this section, we depart from the setting of polynomial dynamics and work in the context of branched coverings.  We introduce model surfaces and model maps, designed to represent restrictions of a polynomial to its basin of infinity.

\subsection{Translation surfaces and horizontal surfaces}

Let $X$ be a Riemann surface, possibly with boundary, and $\omega$ a holomorphic 1-form on $X$.  Away from the zeros of $\omega$, the collection of locally defined functions of the form $\psi(z) = \int_{z_0}^z \omega$ provide a compatible atlas of charts into $\C$.  The ambiguity in the definition of these charts is a complex translation.  It follows that away from the zeros of $\omega$, the length element $|\omega|$ defines a flat Riemannian conformal metric, and this metric extends to a length metric on $X$ with conical singularities at the zeros of $\omega$.  Conversely, given an atlas $\{(U, \psi_U)\}$ on a Riemann surface where the overlap maps between charts differ by translations, the $1$-form $\omega$ on $X$ defined by $\omega = \psi_U^*(dz)$ is globally well-defined.  We call such a pair $(X,\omega)$ a {\em translation surface}.  

A translation surface has natural {\em horizontal and vertical foliations} given by the inverse images of horizontal and vertical lines under the above defined local charts.  These foliations have singularities at the zeros of $\omega$.  At a zero of multiplicity $k$, the metrical has a conical singularity with total angle $2\pi(k+1)$.  

A {\em horizontal translation surface} is a translation surface $(X,\omega)$ for which the overlap maps between charts are translations of the form $z \mapsto z+c$, $c \in \R$.  
On such a surface, there is a globally defined harmonic {\em height function} $G_X: X \to \R$, well-defined up to an additive constant,  given by 
	$$G_X(x) = \int_{x_0}^x \Im \omega.$$  
The connected components of its level sets are the leaves of the horizontal foliation.

\subsection{Model surfaces and local model surfaces.} 
In our applications, we will only encounter horizontal translation surfaces with additional properties.  We single them out as follows.  A {\em model surface} is a connected horizontal translation surface  $(X,\omega, C_X)$ with a distinguished {\em core} $C_X \subset X$ satisfying the following properties:
\begin{itemize}
\item	$X$ is planar (genus 0);
\item the image of the height function $G_X: X\to \R$ is an interval $(a,b)$ with  $\infty \leq a < b \leq \infty$, and for all $c \in (a,b)$ the level sets of the height function $G_X$ are compact and have constant length $\int_{\{G_X = c\}} |\omega| = 2\pi$;
\item the core $C_X$ is of the form $G_X^{-1}[c_0, c_1]$ where $a < c_0 \leq c_1 < b$, and it contains all of the singular leaves of the horizontal foliation on $X$;
\item	 for all $c \geq c_1$, the level set $\{G_X = c\}$ is connected.
\end{itemize}
It follows (see \S 4.3) that $X\setminus C_X$ is a disjoint union of annuli, the {\em outer annulus} $G_X^{-1}(c_1, b)$ and the finitely many {\em inner annuli} with union equal to $G_X^{-1}(a, c_0)$.

A {\em local model surface} is a model surface $(X,\omega, C_X)$ such that $C_X$ consists of a single leaf of the horizontal foliation.

Two model surfaces $(X, \omega, C_X)$ and $(Y, \eta, C_Y)$ are {\em isomorphic} if there is a conformal isomorphism $f: X \to Y$ such that $\omega = f^*(\eta)$ and $f^{-1}(C_Y) = C_X$.  

\begin{figure}
\includegraphics[width=4.5in]{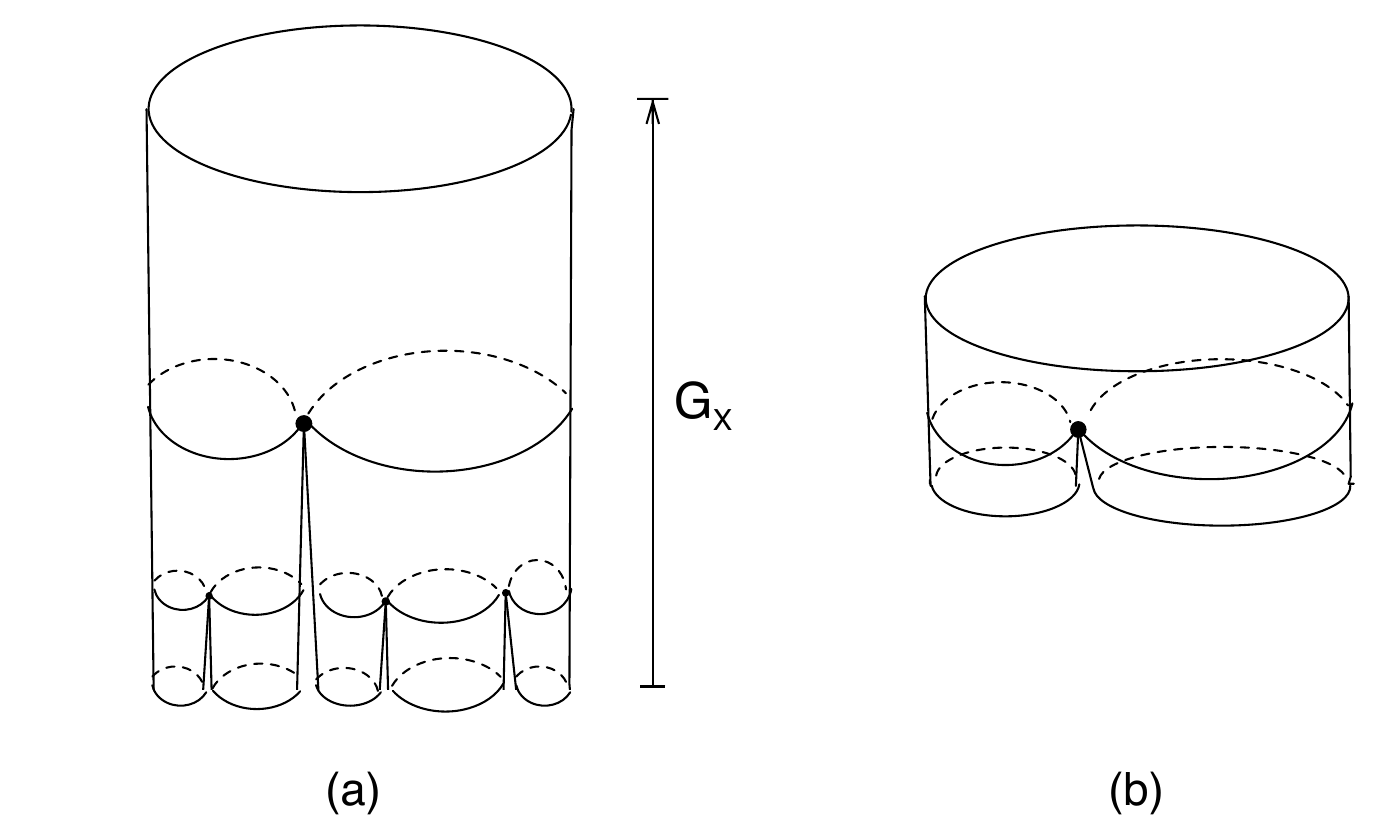}
\caption{(a) A model surface and (b) a local model surface.}
\label{model surface}
\end{figure}

\subsection{Example: local model surfaces without singular leaves}  
\label{subsecn:annuli} Suppose $(X, \omega, L_X)$ is a local model surface, and the core $L_X$ is a single leaf without singularities.  Then $X$ is foliated by horizontal leaves of constant length $2\pi$, so the metric space $(X, |\omega|)$ is a Euclidean cylinder of circumference $2\pi$ and height $h\in (0, \infty]$.  For a suitable choice of positive real constant $c$, the map $\phi: X \to \C$ given by $x \mapsto c \cdot \exp(- i  \int_{x_0}^x \omega)$ defines an isomorphism of horizontal translation surfaces $(X, \omega) \to (Y, i \frac{dz}{z})$ where $Y$ is exactly one of the following subsets of the plane:
\begin{enumerate}
\item  $\{0 < |z| < \infty\}$, 
\item (a) $\{0 < |z| < 1\}$ or (b) $\{1 < |z| < \infty\}$, or 
\item $\{1 < |z| < e^h\}$ if $h<\infty$; 
\end{enumerate}
In case (1), any two choices of core leaf $L_X$ yield isomorphic local model surfaces.   In case (2), let $c$ denote the height of $L_X$.  The distance of the central leaf to the boundary (equivalently, the modulus of the annulus $\{G_X > c\}$ and $\{G_X < c\}$, for cases (a) and (b) respectively) provides a complete invariant for the isomorphism type of $(X, \omega, L_X)$.  In case (3), the modulus of the annulus is $h/2\pi$.  The modulus and the distance from the central leaf to the outer boundary component together provide such an invariant.

\subsection{Example: polynomial pull-back}  \label{poly pullback}
Fix $k \geq 1$, and real numbers $0 < r \leq R < \infty$.  Let $A = \{ z \in \C : r \leq |z| \leq R \}$.  Let $f$ be a polynomial of degree $k$ such that all critical values lie in $A\cup\{0\}$.  Let $X = f^{-1}(\C\setminus\{0\})$, $\omega = i\, k^{-1} f^*(dz/z)$, and $C_X =f^{-1}(A)$.  Then $(X,\omega, C_X)$ is a model surface.  Writing 
	$$f(z) = a\prod_j(z - q_j)^{m_j},$$
we have 
	$$\omega = i\, k^{-1} d(\log f(z)) =  \frac{i}{k} \sum_j \frac{m_j}{z-q_j}  \; dz.$$
In particular, $\omega$ is a meromorphic 1-form with only simple poles; it has residue $i\, m_j/k$ at each finite pole $q_j\in\C$; note $\sum_j m_j/k = 1$.    When $r=R$, $(X, \omega, C_X)$ is a local model surface; the leaf $C_X$ is singular if and only if $f$ has a nonzero critical value.

\subsection{Uniformizing the model surfaces} \label{local model uniformization}
\label{subsecn:uniformization}
Fix a model surface $(X, \omega, C_X)$ and a point $x_0\in X$ lying on a nonsingular vertical leaf.  Suppose the height function $G_X$ has image the interval $(a,b)$, where $ -\infty \leq a < b \leq \infty$.  Delete the vertical leaf through $x_0$, and delete also the set of points in any  singular vertical leaves at and below the heights of zeros of $\omega$.  The resulting subset of $X$  is a connected and simply-connected domain $W$.  The map $\psi: W \to \C$ given by $\psi(x) =  \int_{x_0} ^x \omega$, when post-composed with a real translation, defines an isomorphism from $W$ onto a slit rectangular domain 
$$ R=\{ \theta + i h : 0 < \theta < 2\pi, a < h < b \} - \Sigma $$
where $\Sigma$ is a (possibly empty) finite collection of vertical segments
$$ \Sigma = \bigcup_k \; \, \{\theta_k\} \times \{a < h \leq c_k\} $$
where each $c_k < b$.  
It follows that every model surface can be formed from such a slit rectangular domain $R$, identifying the left side of a slit with real part $\theta_k$ with the right  side of slit with real part $\theta_{\sigma(k)}$ for some permutation $\sigma$ on the set of slits.  (The identification of vertical edges in $R$ must satisfy some obvious planarity conditions which will be treated in more detail in a sequel to this paper.)  Then $\omega = \psi^*(dz)$, so the conformal metric defined by $\omega$ is simply the pullback of the flat Euclidean metric under $\psi$.

\medskip
Next, we describe how a model surface has another Euclidean incarnation, generalizing the polynomial pull-back examples of \S\ref{poly pullback}.

\begin{lemma} \label{pair rep}
Every model surface $(X,\omega, C_X)$ embeds uniquely into a maximal local model surface whose underlying Riemann surface is isomorphic to a finitely punctured plane.  By uniformization, the surface and 1-form is represented by 
	$$\left( \C\setminus\{q_1, \ldots, q_n\} ,\; i\, \sum_j \frac{r_j}{z-q_j} \; dz \right)$$
for some finite set $\{q_1, \ldots, q_n\}$ in $\C$ and real numbers $r_j>0$ such that $\sum r_j = +1$.  Such a representation is unique up to an affine transformation $A \in \Aut \C$.  
\end{lemma}

\proof
The inner and outer annuli of $(X,\omega, C_X)$ can be extended to half-infinite cylinders, so that the height function takes all values in $(-\infty, \infty)$ and the $1$-form extends uniquely to form a new model surface $(\widehat{X}, \omega, C_X)$ which is complete as a metric space.  The ends of $\widehat{X}$ are isomorphic to punctured disks.  The complex structure extends over the punctures to yield a compact Riemann surface $\cl{X}$.  Since $X$ is assumed planar, $\cl{X}$ is homeomorphic to the sphere, hence by the Uniformization Theorem is isomorphic to the Riemann sphere $\Chat$ via an isomorphism $\phi: \Chat \to \cl{X}$ taking infinity to the point at height $+\infty$.  The isomorphism is unique up to precomposition with automorphisms of $\Chat$ that fix $\infty$.

The 1-form $\phi^*\omega=M(z) \, dz$ is holomorphic on the complement in $\C$ of a finite set $\{q_1, \ldots, q_m\}$ of points corresponding to the images of the ends of $\widehat{X}$.   The inner and outer annuli of $\widehat{X}$  are each isometric to a half-infinite Euclidean cylinder of some circumference $r_j$.  The map $\psi_j(x) = \exp(-2\pi i /r_j \int^x \omega)$ then provides an isomorphism (of horizontal translation surfaces) from this annulus to one of the types given in  \S \ref{subsecn:annuli}.  It follows that $M(z) \, dz$ is meromorphic on $\Chat$, the point at infinity is a simple pole with residue $-i$,  and the points  $q_j$ are simple poles with residues $i r_j$ satisfying $\sum_j r_j = 1$.  
\qed

\subsection{External angles.}
Suppose $(X,\omega, C_X)$ is a model surface, and let $e: \widehat{X} \hookrightarrow \C$ be an embedding of the canonical extension given by Lemma \ref{pair rep}. On the outer annulus of $\widehat{X}$, each vertical leaf is nonsingular.  For $z$ large we have, since $\sum_j r_j=1$, that 
\[ e_*(\omega) = i\left(1+O(z^{-1})\right)\frac{dz}{z}.\]
Hence each vertical leaf has a limiting asymptotic argument, $\theta \in \R/2\pi \Z$, at infinity.  Thus, once the embedding $e$ has been chosen, one may speak meaningfully of the {\em external ray of angle $\theta$}  of a model surface $(X,\omega, C_X)$.  These angles coincide with the $\theta$-coordinate of a suitable rectangular representation of $(X, \omega, C_X)$ as in \S\ref{local model uniformization}.

\subsection{Model maps.}
A {\em branched cover} of model surfaces 
	$$f: (Y,\eta, C_Y) \to (X,\omega, C_X)$$
is a holomorphic branched cover $f:Y\to X$ such that 
	$$\eta = \frac{1}{\deg f}\; f^*\omega$$
and $f^{-1}(C_X) = C_Y$.  It follows that any critical values of $f$ must lie in the core $C_X$ and that the outer annulus of $Y$ is mapped by a degree $\deg f$ covering map onto the outer annulus of $X$.

Two model maps $f:(Y,\eta, C_Y)\to (X,\omega,C_X)$ and $g: (Z,\nu, C_Z)\to (X,\omega, C_X)$ over the same base $(X,\omega, C_X)$ are {\em equivalent} if there exists an isomorphism of model surfaces $i$ such that 
$$\xymatrix{ (Y,\eta, C_Y) \ar[dr]_f    \ar[rr]^-{i}& & (Z, \nu, C_Z) \ar[dl]^g \\
			 &(X,\omega, C_X) &   }$$
commutes.    

The proof of the following lemma is a straightforward application of the ideas in the proof of Lemma \ref{pair rep}.  

\begin{lemma} \label{poly extension}
Via the embedding of Lemma \ref{pair rep} applied to both domain and range, every model map is the restriction of a polynomial which is unique up to affine changes of coordinates in domain and range.
\end{lemma}

\medskip
Note in particular that the number of critical values in $X$ of a model map $f: (Y,\eta,C_Y)\to (X,\omega, C_X)$ is at most $(\deg f)-1$.

\subsection{Spaces of model maps.} \label{spaces of local models}
Fix a model surface $(X,\omega, C_X)$ and an embedding into $\C$, as in Lemma \ref{pair rep}; thus $X$ and $C_X$ are regarded as subsets of $\C$.  By Lemma \ref{poly extension}, every model map over $(X,\omega, C_X)$ is the restriction of a polynomial; by precomposing with an automorphism of $\C$, we can assume that the polynomial is monic and centered.  Via this representation, the set of equivalence classes of model maps over $(X,\omega, C_X)$ inherits a topology from the space of monic and centered polynomials.

In detail, let $\MM_k(X,\omega, C_X)$ be the set of equivalence classes of model maps of degree $k$ over $(X,\omega, C_X)$.  Fix an embedding 
	$$e: (X,\omega) \to \left( \C\setminus\{q_1, \ldots, q_n\} , \; i\, \sum_j \frac{r_j}{z-q_j} \; dz \right)$$
as given by Lemma \ref{pair rep}.  
Recall that $e$ is uniquely determined up to postcomposition by an affine transformation.  Let $C_X\subset \C$ denote as well the image of the set $C_X$ under the embedding $e$.  Let $\cP_k(X,\omega, C_X)$ be the collection of monic and centered polynomials of degree $k$ with all critical values contained in the set 
	$$C_X \cup \{q_1, \ldots, q_n\}.$$
Note that the restriction on the location of the critical values implies that the preimage $p^{-1}(e(X))$ is connected for any $p$ in $\cP_k(X,\omega, C_X)$.

\begin{lemma} \label{bijection}
Restriction of polynomials defines a bijection
	$$\cP_k(X,\omega, C_X)/ \< \zeta: \zeta^k=1 \> \to \MM_k(X,\omega, C_X)$$
where the $k$-th roots of unity act on polynomials by precomposition:  $\zeta\cdot p(z) = p(\zeta z)$.
\end{lemma}

\proof
For each polynomial $p$ in $\cP_k(X,\omega, C_X)$, its restriction to the connected subset $p^{-1}(e(X))$ defines a model map 
	$$p: \left( p^{-1}(e(X)), \frac{1}{k}\,  p^*\omega, \, p^{-1}(C_X) \right) \to (X,\omega, C_X)$$
of degree $k$.  Precomposing $p$ by a rotation of order $k$ produces another element of $\cP_k(X,\omega, C_X)$ which is clearly an equivalent local model.  Surjectivity follows from Lemma \ref{poly extension}, and injectivity follows from the uniqueness (up to conformal automorphism) of the extension in Lemma \ref{poly extension}. 
\qed

\medskip
The bijection of Lemma \ref{bijection} induces a topology on $\MM_k(X,\omega, C_X)$, as a quotient space of the subset $\cP_k(X,\omega, C_X)$ of $\cP_k$, the space of all monic and centered polynomials of degree $k$.  While the set $\cP_k(X,\omega, C_X)$ depends on the choice of embedding $e$, the quotient sets $\MM_k(X,\omega, C_X)$ are canonically homeomorphic for any two such choices.  Indeed, suppose $e_1$ and $e_2$ are two embeddings and let $\cP_k^{(1)}(X,\omega, C_X)$ and $\cP_k^{(2)}(X,\omega, C_X)$ be the corresponding sets of polynomials.  The composition $e_2\circ e_1^{-1}$ extends to an affine automorphism $z \mapsto az+b$ of $\C$.  It follows that $p(z) \in \cP_k^{(1)}(X,\omega, C_X)$ if and only if $e_2\circ e_1^{-1} \circ p(a^{-1/k}z) \in \cP_k^{(2)}(X,\omega, C_X)$ for any choice of root $a^{-1/k}$.

\begin{lemma} \label{maximal component}
Fix a model surface $(X,\omega, C_X)$.  The subset 
	$$\MM_k^{k-1}(X,\omega, C_X)\subset \MM_k(X,\omega, C_X),$$
consisting of  model maps represented by polynomials with all $k-1$ critical values in $C_X$, is compact and path-connected.  Any model map in $\MM_k^{k-1}(X, \omega, C_X)$ sends an inner annulus by degree one onto its image.
\end{lemma}

\proof
Let $S$ be the subset of $\cP_k(X,\omega, C_X)$ consisting of polynomials with all $k-1$ critical values in the compact and path-connected set $C_X$.  By Lemma \ref{compactconnected}, $S$ is compact and path-connected.  By Lemma \ref{bijection}, the subset $\MM_k^{k-1}(X,\omega, C_X)$ is homeomorphic to the quotient $S/ \< \zeta: \zeta^k=1 \>$, hence is also compact and path-connected.  
The complement $X\setminus C_X$ is a disjoint union of annuli which neither meet nor surround critical values of the representing polynomial, so the final statement of the Lemma follows.  
\qed

\subsection{Pointed model surfaces and maps}
Let $(X,\omega, C_X)$ be a model surface.  A {\em pointed  model surface} is a quadruple $(X,x,\omega, C_X)$, where $x$ is any point in the outer annulus of $X$.  We consider pointed model maps 
	$$f: (Y,y,\eta,C_Y) \to (X,x,\omega, C_X),$$
i.e. model maps $f: (Y,\eta, C_X) \to (X,\omega, C_X)$ such that $f(y) = x$.  
Two pointed model maps $f, g$ with the same image are {\em equivalent} if there exists an isomorphism $i$ of pointed local model surfaces such that
$$\xymatrix{ (Y,y,\eta, C_Y) \ar[dr]_f    \ar[rr]^-{i}& & (Z, z, \nu, C_Z) \ar[dl]^g \\
			 &(X,x,\omega, C_X) &   }$$
commutes. 
We let $\MM_k(X,x,\omega, C_X)$ denote the set of equivalence classes of these pointed model maps.  As in the non-pointed case, the set can be topologized via an identification with monic and centered polynomials.  Let $\cP_k(X,\omega, C_X)$ be the set of monic and centered polynomials defined in \S\ref{spaces of local models}.  Compare the statement of the following lemma to that of Lemma \ref{bijection}.

\begin{lemma} \label{pointed topology}
Let $(X,x,\omega, C_X)$ be a  pointed model surface.  The canonical projection $\MM_k(X,x,\omega, C_X) \to \MM_k(X,\omega, C_X)$ factors through a bijection $b$ such that the diagram   
	$$\xymatrix{ \MM_k(X,x,\omega, C_X) \ar[d] \ar[r]^-b & \cP_k(X,\omega, C_X) \ar[d] \\
		\MM_k(X,\omega, C_X) & \ar[l]_-r \cP_k(X,\omega, C_X)/ \<\zeta: \zeta^k = 1\>  } $$
commutes, where $r$ is the restriction map of Lemma \ref{bijection}.
\end{lemma}

\proof
Fix an embedding 
$$e: (X,\omega) \to \left( \C\setminus\{q_1, \ldots, q_n\} , \; i\, \sum_j \frac{r_j}{z-q_j} \; dz \right)$$
so that the marked point $x$ lies on a vertical leaf with external angle 0.  For each element $f: (Y,y,\eta, C_Y) \to (X,x,\omega, C_X)$ of $\MM_k(X,x,\omega, C_X)$, choose an extension of the domain so that the marked point $y$ lies on a vertical leaf of external angle 0.  This uniquely determines an element of $\cP_k(X,\omega, C_X)$; denote this element by $b(f)$.  If two pointed local model maps extend to the same polynomial, then they are clearly isomorphic, via an isomorphism which preserves the marked points; this proves injectivity of $b$.  For surjectivity of $b$, note that the restriction of any element $p\in \cP_k(X,\omega, C_X)$ to $p^{-1}(e(X))$ determines an element of $\MM_k(X,x,\omega, C_X)$ with marked point chosen as the unique preimage of $x$ on the external ray of angle 0.  Consequently $b$ is a bijection.  The diagram commutes by construction. 
\qed

\medskip
The bijection $b$ of Lemma \ref{pointed topology} induces a topology on the set $\MM_k(X,x,\omega, C_X)$, making the projection $\MM_k(X,x,\omega, C_X)\to \MM_k(X,\omega, C_X)$ continuous.  The following is then an immediate consequence of Lemma \ref{compactconnected}:

\begin{lemma} \label{maximal pointed component}
Fix a pointed model surface $(X,x,\omega, C_X)$.   The subset $$\MM_k^{k-1}(X,x,\omega, C_X)\subset \MM_k(X,x,\omega, C_X),$$ consisting of pointed model maps with all $k-1$ critical values $C_X$, is compact and path-connected.  Any such model map sends an inner annulus by degree one onto its image.

\end{lemma}

In section \ref{basin models}, pointed models are used to keep track of external angles; a point in the outer annulus of $X$ marks a unique vertical leaf in the foliation of $\omega$.  

\subsection{Gromov-Hausdorff topology on a space of model maps}
Let $(X, x, \omega, C_X)$ be a pointed model surface.  In addition to the algebraic topology inherited as a subset of $\cP_k$ (see Lemma \ref{pointed topology}), the space of model maps $\MM_k(X, x, \omega, C_X)$ also admits a natural Gromov-Hausdorff topology, as for polynomials on their basins of infinity; see \S\ref{basins:topology}.  Two maps $f_j: (Y_j, y_j, \eta_j, C_j) \to (X, x, \omega, C_X)$, $j=1,2$ are called $\epsilon$-close if there is a relation $\Gamma \subset Y_1 \times Y_2$ containing $(y_1, y_2)$ which is nearly surjective and nearly isometric on a neighborhood of the core $C_1\times C_2$ (as in \S\ref{basins:topology}, with the obvious modifications) and with the condition of nearly conjugate replaced with the following.  
If $(y, y') \in \Gamma$, then $\rho_X(f_1(y), f_2(y')) < \epsilon$, where $\rho_X$ is the distance function on $X$ determined by $|\omega|$.  

\begin{lemma}
\label{lemma:two_topologies_on_models} Let $(X, x, \omega, C_X)$ be a pointed model surface.  The Gromov-Hausdorff and algebraic topologies on $\MM_k(X, x, \omega, C_X)$ coincide.
\end{lemma}

\proof The arguments are similar to those used to prove continuity of $\pi: \MP_d \to \cB_d$.  If two polynomials $f_1, f_2 \in \MM_k(X, x, \omega, C_X)$ are close in the algebraic topology, then when uniformized they are close as elements of $\cP_k$.  The topology on $\cP_k$ is by coefficients; or equivalently, locally uniform convergence.  Consequently, for any compact set $K$ in $X\subset \C$, the subsets $f_i^{-1}(K)$ will be close as subsets of the plane.  The $f_i$ and their derivatives are uniformly close on $K$, so the 1-forms $f_i^*(\omega)$ will also be close.  Consequently, $f_1$ and $f_2$ are Gromov-Hausdorff close.  

By compactness of $\MM_k(X, x, \omega, C_X)$ in the algebraic topology, it suffices to observe that $f_1$ and $f_2$ are equivalent as model maps if and only if they are isometrically conjugate on a neighborhood of their cores.  This is clear because the 1-forms on the inner and outer annuli are determined by the core and local degree of the covering.  
\qed

\bigskip\bigskip
\section{Model maps in the basin of infinity}
\label{basin models}

Let $f \in \MP_d$.  In this section we construct model surfaces and model maps in the basin of infinity of $f$.  We define the gluing operation, where a piece of $f|X(f)$ is replaced by a new model map.  We prove the continuity of the gluing operation, and we show that the set $\cS(f,t) \subset \MP_d$ is path-connected for every $f$ and every $t>0$.

\subsection{Forming model surfaces from the basin of infinity}
For any pair of real numbers $0 < a < b < \infty$, a connected component $X$ of $\{a < G_f < b\}$ forms a model surface in the following way.  There are only finitely many singular leaves in $X$.  Let $L$ be the highest singular leaf in $X$ (or choose any leaf if there are no singular leaves). Let $c_1 = G_f(L)$ be the height of the leaf $L$ and let $l$ be the level of $L$, as defined in \S\ref{metric}.   Define
\begin{equation} \label{scaled omega}
	\omega_X = \frac{d^l}{\deg(f^l|L)} \;\omega_f = \frac{2\pi}{\int_L |\omega|} \; \omega_f.
\end{equation}
Let $c_0$ be the smallest height of a singular leaf in $X$ (or equal to the height of $L$ if there are no singular leaves), and set $C_X = \{z\in X: c_0 \leq G_f(z) \leq c_1\}$.  Then the triple $(X,\omega_X, C_X)$ is a model surface.

\subsection{Gluing in new model maps}  \label{gluing}
Fix a polynomial $f\in \MP_d$.  Choose any pair of real numbers $0 < c_0 \leq c_1 < d\cdot c_0$.  Here we specify a collection of pointed model maps which model the restriction of $f$ to $G_f^{-1}[c_0, c_1] \subset X(f)$.  Our goal is to define a process of extracting these models from the basin of $f$ and gluing in new model maps.  

Given $0 < c_0 \leq c_1 < d\cdot c_0$, there exist $0 < a < c_0 \leq c_1 < b$ so that each component $Z_1, \ldots, Z_r$ of the locus $\{a < G_f < b\}$ forms a model surface with core at height $[c_0, c_1]$.  Let $\omega_{Z_i}$  be defined by equation (\ref{scaled omega}) and let $C_{Z_i}$ be the component of $G_f^{-1}[c_0, c_1]$ in $Z_i$; the model surface is $(Z_i, \omega_{Z_i}, C_{Z_i})$.   

\begin{figure}
\includegraphics[width=5in]{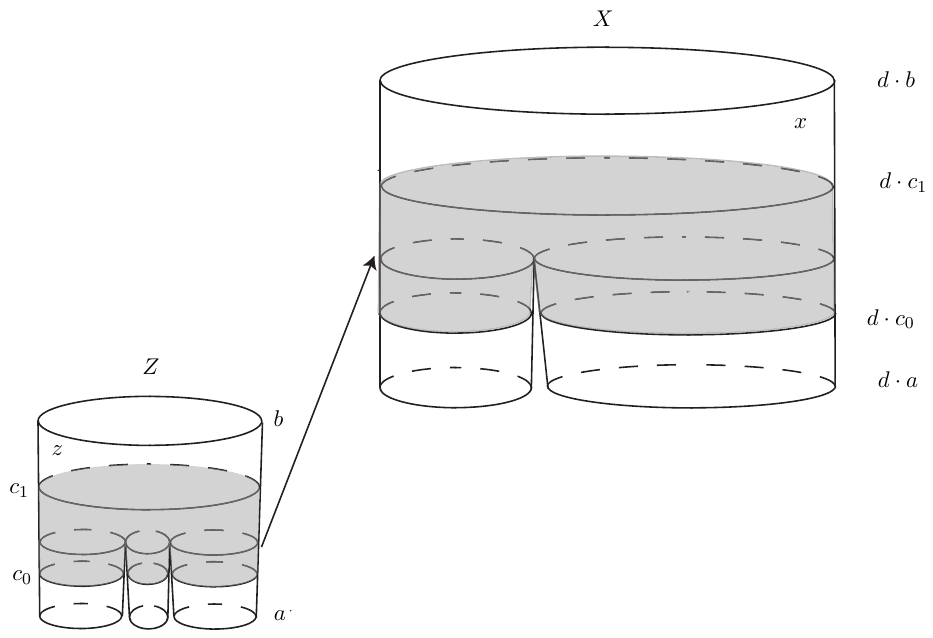}
\caption{Model of $f$ at height $[c_0, c_1]$.}
\label{polynomial model}
\end{figure}

Label the images $X_i = f(Z_i)$; note that we may have $X_i = X_j$ for $i\not=j$.  For each $i$, form the model surface $(X_i, \omega_{X_i}, C_{X_i})$ where $\omega_{X_i}$ is defined by equation (\ref{scaled omega}) and $C_{X_i} = f(C_{Z_i})$.  The restriction of $f$ defines a model map
	$$f|_{Z_i}: (Z_i, \omega_{Z_i}, C_{Z_i}) \to (X_i, \omega_{X_i}, C_{X_i})$$
For each $i$, choose a point $x_i$ in the outer annulus of $X_i$ and let $z_i$ be any preimage of $x_i$ in $Z_i$.  We thus obtain a family of pointed model maps 
	$$f|_{Z_i}: (Z_i, z_i, \omega_{Z_i}, C_{Z_i}) \to (X_i, x_i, \omega_{X_i}, C_{X_i})$$
See Figure \ref{polynomial model}.   Let $k_i = \deg(f|_{Z_i})$.

We now define a sort of inverse procedure, which we call {\em gluing}.  For each $i = 1, \ldots, r$, choose any pointed model map
	$$p_i: (Y_i, y_i, \eta_i, C_{Y_i}) \to (X_i, x_i, \omega_{X_i}, C_{X_i})$$
of degree $k_i$ over the given base $(X_i, x_i, \omega_{X_i}, C_{X_i})$.   The restrictions of $p_i$ and of $f_i$ to the outer annuli of $Y_i$ and of $Z_i$, respectively, are covering maps of degree $k_i$ onto the outer annulus of $X_i$.  There is a unique conformal isomorphism identifying these annuli which sends $y_i$ to $z_i$ and pulls $\omega_{Z_i}$ back to $\eta_i$ on the outer annulus.  Via these identifications, we form a new Riemann surface 
	$$X_a := X_{c_1}(f) \cup \bigcup_{i=1}^r Y_i$$
The 1-form and height function extend to the surface $X_a$ to yield a one-form $\tilde{\omega}$ and a height function $\tilde{G}: X_a \to \R$ with image equal to $(a, \infty)$.  The map $f|X_{c_1}(f)$ extends holomorphically to a new self-map 
	$$\tilde{f}: X_a \to X_a$$
which agrees with $p_i$ on $Y_i$ and satisfies $\tilde{G}(\tilde{f}(z)) = d\cdot\tilde{G}(z)$.  We say that $\tilde{f}$ is obtained from $f$ by {\em gluing in the model maps $p_i$ at height $[c_0, c_1]$}.  

Note that different choices of the constants $a, b$ (chosen in the second paragraph of this section) yield, after suitable extensions, isometrically conjugate maps.  If e.g. $a' < a$, the locus $\{ a' < G_f < a\}$ consists of annuli.  The map $\tilde{f}: X_a \to X_a$ is affine near the lower boundary, and therefore extends canonically over these annuli to yield the map $\tilde{f}': X_{a'} \to X_{a'}$.

\subsection{Continuity of gluing} \label{glue continuity}
In the gluing construction of \S\ref{gluing}, we are particularly interested in the case where each $p_i$ has exactly $k_i-1$ critical values in the core $C_{X_i}$, counted with multiplicity. We show that $\tilde{f}$ extends uniquely to a polynomial in the shift locus of $\cM_d$.  Recall the notation $\pi_t: \MP_d \to \cB_{d,t}$ from \S\ref{basins:topology}; in the space $\cB_{d,t}$, two polynomials $f_i, i= 1,2$ are equivalent if their restrictions to $X_t(f_i) = \{G_{f_i} >t\}$ are conformally conjugate.

\begin{prop}
\label{prop:gluing yields polynomial}
For each $f \in \MP_d$ and any pair of real numbers $0 < c_0 \leq c_1 < d\cdot c_0$, let  
	$$(f_i) \in \prod_i \MM_{k_i}(X_i, x_i, \omega_{X_i}, C_{X_i})$$
be a pointed model representation of the restriction $f| G_f^{-1}[c_0, c_1]$ as constructed above.  Then gluing at height $[c_0, c_1]$ defines a continuous map
\[\glue: 
	\prod_i \MM_{k_i}^{k_i-1}(X_i, x_i, \omega_i, C_{X_i}) \to \cS_d\cap \{g : m(g) \geq c_0 \}\]
such that $\pi_{c_1} \circ \glue(p_1, \ldots, p_r) = \pi_{c_1}(f)$ for any choice of $p_1, \ldots, p_r$.  

If $f \in \cS_d$ and $m(f)> c_1$, then the domain of $\glue$ is the single point $(f_1, \ldots, f_r)$, and $\glue(f_1, \ldots, f_r)=f$.  That is, $f$ is determined by its restriction $f|X_{c_1}(f)$.    
\end{prop}

\proof  We adopt the notation as in the definition of gluing in \S\ref{gluing}.  We first claim that gluing defines a polynomial in the shift locus.  Choose model maps $(p_i)$ to glue to $f$ at height $[c_0, c_1]$.  Let $\tilde{f}: X_a\to X_a$ be the extended self-map of the surface $X_a$ as in the definition of gluing.  We argue inductively that $\tilde{f}$ extends to a proper, holomorphic self-map of a rigid planar Riemann surface to itself.  

Let $\{W_j\}$ denote the collection of inner annuli for the model surfaces $(Y_i, y_i, \eta_i, C_{Y_i})$, so the $\{W_j\}$ are the connected components in $X_a$ of $\{a < \tilde{G} < c_0\}$.  Let $V_j = \tilde{f}(W_j)$ denote the images of $W_j$.  Thus, each $V_j$ is an annulus whose height function has image $(d\cdot a, d\cdot c_0)$.  By our choices of $p_i\in \MM_{k_i}^{k_i-1} (X_i, x_i, \omega_{X_i}, C_{X_i})$, all of the inner annuli $W_j$ map with degree 1 to their images.  Each $V_j$ is an annulus with height in the interval $(da, d\cdot c_0)$.  

For each $j$, construct a model surface with outer annulus $V_j$.  Namely, we may take the connected components of $\{a < \tilde{G} < d c_0\}$ in $X_a$ as our model surface $(V_j', \omega_{V'_j}, C_{V'_j})$, with core at height $[c_0, da]$.  We now redo the gluing procedure with this collection of model surfaces as the base.  Indeed, there is a {\em unique} model surface with outer annulus $W_j$ which is isomorphic to $(V_j', \omega_{V'_j}, C_{V'_j})$, and the degree 1 restriction $\tilde{f}: W_j \to V_j$ extends uniquely to the new model surface.  In this way, we extend $\tilde{f}$ holomorphically to a new surface
	$$\tilde{f} : X_{a/d} \to X_{a/d}$$
with a height function $\tilde{G}: X_{a/d} \to (a/d, \infty)$ satisfying $\tilde{G}(\tilde{f}(z)) = d\cdot\tilde{G}(z)$.  We now repeat the extension procedure by setting the new annular components $W_j$ to be the connected components of $\{a/d < \tilde{G} < c_0/d\}$.  By induction, $\tilde{f}$ extends to a proper, degree $d$, holomorphic self-map 
	$$\tilde{f}: X \to X$$ 
of a planar Riemann surface $X$ to itself.  

At every step of the induction, the annuli $\{W_j\}$ map by degree 1 to their images, so we see easily that the Riemann surface $X$ satisfies the infinite-modulus condition (see \S\ref{basins:shift}).  It is therefore rigid, and there exists a conformal embedding $X\hookrightarrow\C$, sending $\infty$ to $\infty$, unique up to postcomposition by an affine transformation.  We may conclude that $\tilde{f}$ extends to a polynomial $g:\C\to\C$, unique up to affine conjugation.  We set 
	$$g=\glue(p_1, \ldots, p_r).$$  
By construction, every critical point of $g$ lies in the basin of infinity $X(g)$, so $g \in \cS_d$.  The height function $\tilde{G}$ coincides with the escape-rate function $G_g$ on $X(g)$, and the restrictions $g|X_{c_1}(g)$ and $f|X_{c_1}(f)$ are holomorphically conjugate.   As each $p_i$ has all $k_i-1$ critical points in the core of $X_i$, we may also conclude that $G_g(c) \geq c_0$ for all critical points $c$ of $g$.  

If $m(f)>c_1$, then there is a unique choice for $(p_i)$ at the first stage of gluing, as each $p_i$ defines an isomorphism to its image.   The uniqueness of the extension to a polynomial implies that the restriction $f|X_{c_1}(f)$ determines the conformal conjugacy class of $f$.

We now prove continuity of the map $\glue$.  Fix $\eps, s >0$ such that $s<c_0$, and fix $g$ in the image of $\glue$, so $g\in \cS_d$ with $m(g) \geq c_0$.  We aim to show that the preimage of the Gromov-Hausdorff neighborhood $U_{s,\eps}(g)$ is open in $\prod_i \MM_{k_i}^{k_i-1} (X_i, x_i, \omega_{X_i}, C_{X_i})$.  By Lemma \ref{lemma:two_topologies_on_models}, we may work with the Gromov-Hausdorff topology on the space of model maps.

Let $p = (p_1, \ldots, p_r) \in \prod_i \MM_{k_i}^{k_i-1} (X_i, x_i, \omega_{X_i}, C_{X_i})$ be any point sent to $g$.  Let $V_\eps = \prod_i V_\eps(p_i)$ be the Gromov-Hausdorff neighborhood of $(p_1, \ldots, p_r)$ consisting of all model maps that are $\eps$-close to $p_i$ for each $i$.  From the definition of the topology, for any $q \in V_\eps$, there is a relation $\Gamma_q$ between the domain of $q$ and that of $p$ which shows the model maps are $\eps$-close.  Let $\Delta$ be the identity relation (the diagonal) in $X_{c_1}(f) \times X_{c_1}(f)$.    Recalling that $a$ denotes the minimal height of the domains of the model maps $q$, we see that then $\Gamma_q$ and $\Delta$ together form a relation in $X_a(\glue(q))\times X_a(g)$ which shows that
	$$\glue(q) \in U_{a, \eps}(g)$$
for all $q\in V_\eps$. 

In fact, we will see that $V_\eps$ is sent by $\glue$ into the neighborhood $U_{s,\eps}(g)$ for any $s < a$.  The extension to level $\{\tilde{G}>a/d\}$ is uniquely determined, by gluing in degree 1 model maps.  The relation $\Gamma_q$ between the domain of $q$ and that of $p$ determines a relation on the domains of the model maps at this lower height; the distance $\eps$ is the same.  As glued maps, the distance between the glued images can only decrease, as the metric on a basin of infinity is scaled by $1/d$ with every preimage.  Continuing inductively, we see that
	$$\glue(q) \in U_{a/d^n, \, \eps}(g)$$
for all positive integers $n$ and all $q\in V_\eps$.  We conclude that the map $\glue$ is continuous.\qed

\subsection{Consequences of the gluing construction}
Proposition \ref{prop:gluing yields polynomial} implies the following facts about the projections $\cS_d \to \cB_{d,t}$; precise statements appear in the following theorem.  Suppose $f \in \cS_d$.   If $t>0$ is sufficiently small, then the restriction of $f$ to $X_t(f)$ determines $f$ uniquely.  Second, any map $g$ satisfying $m(g)\geq t$ and such that $g|X_t(g)$ is holomorphically conjugate to $f|X_t(f)$ is obtained by such gluings, and the totality of such maps is connected.  Lastly, as long as the combinatorial data (the number of components of $\{G_f = t\}$) remains constant, gluings can be transported continuously along one-parameter families $f_s$.

\begin{theorem}
\label{thm:proved_by_gluing} 
Let $t>0$.  
\begin{enumerate}
\item The restriction $\pi_t| \cS_d \cap \{g: m(g)>t\} \to \cB_{d,t}$ is a homeomorphism onto its image.   
\item The restriction $\pi_t| \cS_d \cap \{g: m(g) \geq t\} \to \cB_{d,t}$ is surjective, and the fibers are path-connected.   
\item Suppose $f_s, s \in [0,1]$ is a continuous path in $\MP_d$, and $t>0$ has the following stability property: there exist $r \in \N$ and real numbers $a< t < b$ such that for all $s \in [0,1]$, the locus $\{ a < G_{f_s} < b\}$  consists of $r$ annular components.  Suppose $g_0 \in \cS_d \cap \{g: m(g) \geq t\}$ is given, and $\pi_t(g_0)=\pi_t(f_0)$.  Then there exists a continuous path $s \mapsto g_s \in  \cS_d \cap \{g: m(g)\geq t\}$ starting at $g_0$ such that $\pi_t (g_s)=\pi_t(f_s)$ for all $s \in [0,1]$.
\end{enumerate}
\end{theorem}

\noindent
The stability hypothesis in (3) implies that the indicated components, in the associated conformal metrics,  form a family of Euclidean annuli whose isometry types are constant as $s$ varies.  

Given $f$ in the shift locus, we define 
	$$\cS(f,t) = \cS_d \cap \{ g: m(g) \geq t\} \cap \pi_t^{-1}(\pi_t(f)),$$ 
the fiber of the restriction in (2) containing $f$.  In words, the set $\cS(f,t)$ consists of all maps $g$ holomorphically conjugate to $f$ above height $t$ and satisfying $m(g) \geq t$.  The structure of $\cS(f,t)$ will play a role in the proof of Theorem \ref{pi}.  For later reference, we state the following corollary explicitly:

\begin{cor}  \label{S is path-connected}
For any $f$ and any $t>0$, the set $\cS(f,t)$ is path-connected.
\end{cor}

Another important immediate consequence is:

\begin{prop}
\label{prop:shift locus is dense}
The shift locus is dense in $\cB_d$. 
\end{prop}

\proof
Suppose the polynomial $f$ represents an element of $\cB_d$,  let $\epsilon>0$ and suppose $t>0$ satisfies $1/t > M(f)$.  By Theorem \ref{thm:proved_by_gluing}(2), there is a polynomial $g \in \cS_d$ for which $g|X_t(g)$ is holomorphically, hence isometrically, conjugate to $f|X_t(f)$.  Thus $g \in U_{t, \epsilon}(f)$.\qed

\subsection{Proof of Theorem \ref{thm:proved_by_gluing}}
For the proof, it will be more convenient to work with the space $\cP_d$ of monic and centered polynomials, so that  each basin of infinity has well-defined external angles.  In particular, any $f \in \cP_d$ fixes exactly $d-1$ distinct external rays.   Let $\tilde\cB_d$ be the set of conformal conjugacy classes of monic, centered polynomials restricted to their basins of infinity, where now the conjugacy is required to have derivative $1$ at infinity.    As a set, $\tilde\cB_d $ is the set of equivalence class of triples $(f, X(f), \theta_f)$, where $\theta_f$ is one of the $d-1$ external rays  that are fixed under $f$, and where two triples $(f, X(f), \theta_f), (g, X(g), \theta_g)$ are equivalent if there is a holomorphic conjugacy from $f$ on $X(f)$ to $g$ on $X(g)$ sending $\theta_f$ to $\theta_g$.  We equip $\tilde\cB_d $ with the smallest topology such that the natural projection $\tilde\cB_d  \to \cB_d$ is continuous.  More concretely: an $\epsilon$-conjugacy in this setting has the same definition as for the Gromov-Hausdorff topology in \S\ref{basins:topology}, with the following additional requirement.  Observe that if $1/t > M(f)$ then the set $\{G_f= 1/t\} \cap \theta_f$ is a singleton $x_f$.   We require that an $\epsilon$-conjugacy $\Gamma$ send $x_f$ to $x_g$, i.e. $(x_f, x_g) \in \Gamma$.   We refer to this topology as the Gromov-Hausdorff topology on $\tilde\cB_d $.  
The arguments showing that the projection $\pi: \MP_d \to \cB_d$ is a homeomorphism on the shift locus (Lemma \ref{homeo}) immediately show that the projection $\tilde{\pi}: \cP_d \to \tilde\cB_d$ is also a homeomorphism on the corresponding shift locus $\tilde\cS_d \subset \cP_d$ equipped with its algebraic topology inherited from the polynomial coefficients.  Finally, we define $\tilde\cB_{d,t}$ analogously.

Given an element $f$ in $\cP_d$, consider its restriction to $X_t(f)$. Define $\tilde \cS(f,t)$ in $\cP_d$ to be the set of polynomials $g\in \cP_d$ with $g|X_t(g)$ conjugate to $f|X_t(f)$ via a conformal isomorphism with derivative 1 at infinity, and such that $m(g)\geq t$.  Then, for each polynomial $g\in \tilde \cS(f,t)$, there is a {\em unique} isomorphism $\phi_g: X_t(g) \to X_t(f)$ conjugating $g$ to $f$ and sending the ray of angle 0 for $g$ to that of $f$.  

We now establish (1).  The indicated restriction is continuous by Lemma \ref{continuous}, and injective by Proposition \ref{prop:gluing yields polynomial}.  By invariance of domain, the conclusion follows. 

Next, we prove (2).   The surjectivity conclusion follows immediately from Proposition \ref{prop:gluing yields polynomial}.   The path-connectivity of the fibers in (2) will follow once we establish that the corresponding set  $\tilde \cS(f,t)$ in $\cP_d$ is path-connected.  This is what we prove below.

Fix $f\in \cP_d$ and let $t>0$ be arbitrary.  Suppose $\{G_f=t\}$ has $r$ components, and let $X_i^f, Z_i^f$ be local model surfaces as constructed in \S\ref{gluing}, where $c_0 = c_1 = t$.  That is, the leaves $\{G_f = t\}$ form the core $L_{X_i}$ of the $Z_i^f$.  For each $i=1, \ldots, r$, fix a choice of points $x_i$ in the outer  annuli of the local model surfaces $X^f_i$, and fix a choice $z_i$ of their preimages under $f$ in the outer annuli of the local model surfaces $Z^f_i$.   

Suppose $g \in \tilde \cS_d$, $m(g)\geq t$ and $\pi_t(g)=\pi_t(f)$, so that $g$ lies in the fiber over $(X_t(f), f, |\omega_f|) \in \cB_{d,t}$.    Let $\phi_g: X_t(g) \to X_t(f)$ be the unique holomorphic conjugacy as in the discussion above.  For each $i$, let $Z_i^g$ be the component of $\{a < G_g < b\}$ whose outer annulus is the image under $\phi_g^{-1}$ of the outer annulus of $Z_i^f$, let  $z_i^g = \phi_g^{-1}(z_i)$, let $\eta_i^g = \phi_g^*(\omega_i)$, and let $L_{Z_i}^g$ be the component of $\{G_g = t\}$ contained in $Z_i^g$.  Then the restriction $g_i$ of $g$ to $Z_i$ followed by the isomorphism $\phi_g$ yields a pointed model map 
\[ \phi_g \circ g_i: (Z_i^g, z_i^g, \eta_i^g, L_{Z_i}^g) \to (X_i, x_i, \omega_i, L_{X_i});\]
since $g \in \tilde \cS(f,t)$, all inner annuli of $Z_i$ map under $g_i$ by degree one.  
We obtain in this way a well-defined map 
\[ \cL: \tilde \cS(f,t) \to \prod_i \MM_{k_i}^{k_i-1}(X_i, x_i, \omega_i, L_{X_i}).\]
The right-hand side is compact and path-connected by Lemma \ref{maximal pointed component}.  The remainder of the proof is devoted to establishing that $\cL$ is in fact a homeomorphism.

Observe that the proof of Proposition \ref{prop:gluing yields polynomial} can be adapted so that the continuous map $\glue$ is taking values in $\cP_d$, the space of monic and centered polynomials.  Indeed, if we begin with $f\in \cP_d$ with its distinguished external angle $\theta = 0$, we require that the glued and extended map $\tilde{f}: X\to X$ be embedded into $\C$ so that the distinguished vertical leaf is sent to angle $\theta=0$.  In this way Proposition \ref{prop:gluing yields polynomial} yields a continuous map 
  	$$\glue:  \prod_i \MM_{k_i}^{k_i-1}(X_i, x_i, \omega_i, L_{X_i}) \to \tilde \cS(f,t)  $$
so that $\cL\circ \glue$ is the identity.  Hence $\cL$ is surjective.  By compactness of the domain of $\glue$, it suffices to prove that $\cL$ is injective.

Suppose $\pi_t(g_1)=\pi_t(g_2)=\pi_t(f)$, $m(g_1)\geq t$, $m(g_2) \geq t$, and $\cL(g_1) = \cL(g_2)=(p_1, \ldots, p_r)$.  Let $\phi_{g_i}: X_t(g_i) \to X_t(f)$ be the isomorphisms defined above.  The isomorphisms yielding equality of pointed local model maps implied by the condition $\cL(g_1)=\cL(g_2)$  glue to the isomorphism $\phi_{g_2}^{-1}\circ \phi_{g_1}$ from $X_t(g_1)$ to $X_t(g_2)$; therefore above some height $a$ with $a<t$ the polynomials $g_1, g_2$ are conjugate  via a conformal isomorphism with derivative 1 at infinity.  By part (1) of Theorem \ref{thm:proved_by_gluing}, the polynomials $g_1$ and $g_2$ are then affine conjugate on $\C$; by construction, there is an isomorphism with derivative 1 at infinity.  Thus $g_1 = g_2 \in \cP_d$ and the proof of (2) is complete.  
\medskip

\subsection*{Remark.}  An alternative, more intrinsic  proof of Theorem \ref{thm:proved_by_gluing}(2) may be given along the following lines, using Lemma \ref{lemma:two_topologies_on_models}.    The surgery constructions in \cite[Section 8]{Eskin:Masur:Zorich} are affinely natural.  This shows that branch values of local model maps in $L_{X_i}$ can be continuously pushed through zeros of $\omega_i$ in $X_i$, and that these branch values can also be so pushed so as to coalesce together to a single branch value, as in the proof of Lemma \ref{compactconnected}.  
\medskip

We now prove the path-lifting claim (3) of Theorem \ref{thm:proved_by_gluing}.  We will derive the conclusion by proving the corresponding statement for monic centered polynomials.  Let $\pi_t: \cP_d \to \tilde \cB_{d,t}$ denote the corresponding projection.  
The path $f_s$ lifts to a path $f_s^\sharp$ in $\cP_d$.  We will show the existence of a lift of a path $g_s^\sharp \in \tilde \cS_d$ for which $\pi_t(f_s^\sharp) = \pi_t(g_s^\sharp)$, starting from an arbitrary given point $g_0^\sharp \in \tilde \cS(f_0, t)$.
To avoid burdensome notation, we now drop the sharp symbols $\sharp$; thus $f_s, g_s$ denote elements of $\cP_d$ and $\tilde \cS_d$, respectively.
 
In this paragraph, we extract from $f_s$ a continuous family of data for the definition of gluing.  
The stability hypothesis implies that there exists $a<t<b$ such that the components $X_i^{f_s}$ of $\{d\cdot a < G_{f_s} < d\cdot b\}$ and $Z_i^{f_s}$ of $\{ a < G_{f_s} < b\}$ comprise a family of $r$ annuli of constant isometric type as $s$ varies,  that the degrees $k_i$ by which the outer annulus of $Z_i^{f_s}$ maps to that of $X_i^{f_s}$ are also constant.   For each $i$ let $L_{X_i}^{f_s}$ be the leaf $\{G_{f_s}=d\cdot t\} \cap X_i^{f_s}$ and similarly define $L_{Z_i}^{f_s} = \{G_{f_s} = t\} \cap Z_i^{f_s}$.  
In the remainder of this paragraph, we show how to continuously choose the points $x^s_i$ and $z^s_i$ needed to define gluing.   By compactness, there exists $M_0>0$ for which $M(f_s) < M_0$ for all $s$.   Since the $f_s$ are monic and centered, for each $s$, there is unique holomorphic conjugacy $\phi_s: \{G_{f_s} > M_0\} \to \{G_{f_0}> M_0\}$ conjugating $f_s$ to $f_0$ and tangent to the identity at infinity.  Choose an integer $l$ so that $d^l \cdot t > M_0$, and set $A_0 = \{d^l\cdot t < G_{f_0} < d^l \cdot b\}$ .  Let $A_s = \phi_s^{-1}(A_0)$, so that $A_s = \{d^l \cdot t < G_{f_s} < d^l\cdot b\}$.  Choose arbitrarily $x' \in A_0$.  Note that $\{t<G_{f_s} < b\}\cap X_i^{f_s}$ is the outer annulus of $X_i^{f_s}$.  The hypothesis on the path $f_s$ and the height $t$ implies that for each $i$, the map  $\id \times (\phi_s \circ f_s^{l-1}): [0,1] \times X_i^{f_s} \to [0,1] \times A_0$ is an unramified covering.  It follows that for each $i$ we may choose a continuous family $x_i^{f_s}$ of preimages of $x'$ under this covering.  By similar reasoning, we may choose preimages $z_i^{f_s}$ of $x_i^{f_s}$ under $f_s$ in $Z_i^{f_s}$ continuously.  Defining the one-forms $\eta_i^{f_s}$ as in the definition of gluing,  we have constructed for each $s \in [0,1]$ and each $i$ a local model map 
\[ f_{i,s}: (Z_i^{f_s}, z_i^{f_s}, \eta_i^{f_s}, L_{Z_i}^{f_s}) \to (X_i^{f_s}, x_i^{f_s}, \omega_i^{f_s}, L_{X_i}^{f_s}).\]
The stability assumption on $t$ implies that for each $s$ and each $i$, the central leaves $L_{X_i}^s$ are nonsingular.  Thus, the isomorphism types of pointed local model base (image) spaces $(X_i^{f_s}, x_i^{f_s}, \omega_i^{f_s}, L_{X_i}^{f_s})$ are independent of $s$, and so for each $i$, the spaces $\MM_{k_i}^{k_i-1}(X_i^{f_s}, x_i^{f_s}, \omega_i^{f_s}, L_{X_i}^{f_s})$ is canonically identified with $\MM_{k_i}^{k_i-1}(X_i^{f_0}, x_i^{f_0}, \omega_i^{f_0}, L_{X_i}^{f_0})$.  

From the proof of (2) above, gluing gives a  homeomorphism 
\[ \tilde \cS(f_s, t) \to \prod_i \MM_{k_i}^{k_i-1}(X_i^{f_s}, x_i^{f_s}, \omega_i^{f_s}, L_{X_i}^{f_s});\]
composing with the isomorphism of the last paragraph, we have that for each $s$, we have a homeomorphism to a fixed space 
\[ \tilde \cS(f_s, t) \to \prod_i \MM_{k_i}^{k_i-1}(X_i^{f_0}, x_i^{f_0}, \omega_i^{f_0}, L_{X_i}^{f_0}).\]
For $s \in [0,1]$ let  $\cL_{1,s}$ be the inverse of this homeomorphism.  
By assumption, $g_0 \in \tilde \cS(f_0, t)$.  Let $(p_1, \ldots, p_r) = \cL(g_0)$ be the image of $g_0$ under the homeomorphism $\cS(f_0,t) \to \prod_i \MM_{k_i}^{k_i-1}(X_i^{f_0}, x_i^{f_0}, \omega_i^{f_0}, L_{X_i}^{f_0})$.  Finally, set $g_s = \cL_{1,s}(p_1, \ldots, p_r)$.   By construction, $\pi_t(g_s)=\pi_t(f_s)$.  It is clear that $g_s$ varies continuously in the Gromov-Hausdorff topology, since the restrictions $\pi_t(f_s)$ vary continuously, the local models $(p_1, \ldots, p_r)$ that are glued to the $f_s$ are constant, and points in the definition of pointed local models that define the gluing vary continuously.  

This concludes the proof of Theorem \ref{thm:proved_by_gluing}.\qed

\bigskip\bigskip\section{Proof that $\pi$ is monotone}

In this section, we conclude the proof of Theorem \ref{pi}.  It remains to show that the projection 
	$$\pi: \MP_d \to \cB_d$$
is monotone.  Recall this means that its fibers are connected.  

\subsection{The set $\cB(f,t)$.} \label{basins:BS}
Fix a polynomial $f$ and a positive real number $t$.  Recall that 
	$$M(f) = \max\{G_f(c) : f'(c)=0\},  \qquad   m(f) = \min\{G_f(c): f'(c) = 0\},$$
and 
	$$X_t(f) = \{z\in X(f):  G_f(z) > t\}.$$	
We defined 
	$$\pi_t: \cM_d \to \cB_{d,t}$$
to be the projection to the space of conformal conjugacy classes of the restriction $f| X_t(f)$.  Let $\cB(f,t)$ be the collection of all polynomials $g$ in $\MP_d$ such that $g | X_t(g)$ is conformally conjugate to $f|X_t(f)$, i.e. $\cB(f,t)$ is the fiber of $\pi_t$ containing $f$.  Recall that we have defined $\cS(f,t)$ to be the set of all $g\in \cB(f,t)$ for which the minimal critical escape rate satisfies $m(g) \geq t$.  

If $t$ is large enough so that $t\geq M(f)$, then (cf. \S 1.3) $\cB(f,t)=\cB(t)$ consists of all polynomials $g$ with $t \geq M(g)$ and is known to be a closed cell.  In degree 2, for the family $z^2+c$, we have the following dichotomy:  
\begin{itemize}
\item	If $t < M(f)$, then $\cB(f,t) = \cS(f,t) = \{f\}$.
\item	If $t \geq M(f)$, then $\cS(f,t)$ is the equipotential curve $\{c: G_c(0)=t\}$ around the Mandelbrot set,  and $\cB(f,t)$ is the closed ball it bounds.  
\end{itemize}
In every degree, when $f$ is in the shift locus and $t$ is small enough that $t<G(c)$ for all critical points $c$ of $f$, then $\cB(f,t) = \cS(f,t) = \{f\}$.

\subsection{Deforming the basin of infinity}
In the next lemma, we use a ``pushing deformation" to show that $\cB(f,t) \cap \cS_d$ is connected.  The construction is similar to the pushing deformation of \cite[\S 4.2]{Blanchard:Devaney:Keen}; in their case, they push critical values down to smaller heights, while we push critical values up along external rays.  Certain deformations require a change in the local topology of the translation structure, like moving through a stratum of $\cB_d$ defined by prescribing the multiplicities of zeros of the 1-form $\omega_f$; compare \cite[Section 8]{Eskin:Masur:Zorich}.

\begin{lemma} \label{paths}
For any $f\in \mathcal{S}_d$ and any $t>0$, there is a path contained in $\cB(f,t)$ joining $f$ to a point in $\cS(f,t)$.  Furthermore, such a path exists with the following properties: (i) the path may be parameterized as $h \mapsto f_h$, where $m(f) \leq h \leq t$, (ii) $f_{m(f)} = f$ and $f_t \in \cS(f,t)$, and (iii) $f_h \in \cS(f, h)$ for all $m(f) \leq h \leq t$.  
\end{lemma}

\proof If all critical points of $f$ have height at least $t$ then already $f \in \cS(f,t)$.  So suppose 
 $f$ has critical points below height $t$, so that $m(f) < t$ is the height of the lowest critical point.  We will ``push" the lowest critical values from the level curves $\{G_f = d\cdot m(f)\}$ up along their external rays in a continuous fashion, without changing the restriction $f|X_t(f)$, until all critical values have height $\geq d\cdot t$.  
 
Choose a finite sequence of heights $m(f) = h_0 < h_1 < \cdots < h_N = t$ so that $h_{j+1}/h_j < d$ for all $j$.  For each $j$ we will glue in a continuous family of model maps to glue in to $f$ at height $[h_j, h_{j+1}]$.   

Beginning with $j=0$, choose pointed model maps 
	$$f_i: (Z_i, z_i, \eta_i, C_{Z_i}) \to (X_i, x_i, \omega_i, C_{X_i})$$ 
for $f$ with core of $Z_i$ at height $[h_0, h_1]$ as in \S\ref{gluing}.  For each critical value $v$ in $X_i$, there is at least one (and possibly several) vertical leaf segment $\lambda_v$ containing $v$, parameterized by height in $[h_0, h_1]$.  We aim to construct a path of pointed model maps 
	$$p_i^h \in \MM_{k_i}^{k_i-1} (X_i, x_i, \omega_i, C_{X_i}) $$
so that $p_i^{h_0} = f_i$ and the critical values of $p_i^h$ lie on $\lambda_v$ at heights $\geq d\cdot h$.  In particular, via the gluing map of Proposition \ref{prop:gluing yields polynomial}, the polynomial 
	$$f_h = \glue(p_1^h, \ldots, p_r^h)$$
will lie in $\cS(f,h)$ for all $h\in [h_0, h_1]$.  We repeat the process for each $j = 0, \ldots, N-1$ to complete the proof of the Lemma.  

Indeed, recall that $\MM_{k_i}^{k_i-1} (X_i, x_i, \omega_i, C_{X_i})$ can be identified with a subset of the monic and centered polynomials $\cP_{k_i}$ (Lemma \ref{pointed topology}).  Further, by Lemma \ref{cvmap}, the map from polynomials in $\cP_{k_i}$ to their collections of critical values has the path-lifting property.  Therefore, we may begin with the path of critical values satisfying the height conditions we desire, each staying on its vertical leaf segment $\lambda_v$, and we may lift it to a path in $\MM_{k_i}^{k_i-1} (X_i, x_i, \omega_i, C_{X_i}) \subset \cP_{k_i}$.  This produces the desired paths $p_i^h$.  The results about gluing in Proposition \ref{prop:gluing yields polynomial} guarantee that the resulting polynomial $f_h$ is in $\cS(f,h)$ for all $h$.  
\qed

\medskip\noindent{\bf Remark.}  It can be seen from the proof of Lemma \ref{paths} that the ``pushing-up" deformation is {\em canonical} unless the moving critical values encounter zeros of $\omega$.  That is, the path is uniquely determined except when the lowest critical values are pushed up through critical points of $f$ or any of their iterated preimages.  Note, however, that if a choice is made at height $t_0 < t$, the path-connectedness of $\cS(f,h)$ by Corollary \ref{S is path-connected} implies that different choices can themselves be connected by paths within $\cB(f,t) \cap \cS_d$.

\begin{cor} \label{shift B connected}
For any $f$ in $\MP_d$ and $t>0$, the intersection of $\cB(f,t)$ with the shift locus $\cS_d$ is path-connected.  In particular, the shift locus is connected.
\end{cor}

\proof
Fix $f$.  It follows immediately from the definition that $\cB(g,t) = \cB(f,t)$ if and only if $g\in \cB(f,t)$.  Similarly, $\cS(f,t)=\cS(g,t)$ if and only if $g\in \cB(f,t)$.  Thus, we may choose any element $g\in \cB(f,t) \cap \cS_d$ and apply Lemma \ref{paths} to find a path from $g$ to $\cS(g,t)=\cS(f,t)$ contained in $\cB(g,t)=\cB(f,t)$.  As $\cS(f,t)$ is path-connected by Corollary \ref{S is path-connected}, we conclude $\cB(f,t) \cap \cS_d$ is path-connected.  Since the shift locus is an increasing union of sets of the form $\cB(f,t) \cap \cS_d$ where $t>0$ and $1/t>M(f)$, the shift locus is connected.  
\qed

\medskip
Recall that a Gromov-Hausdorff basis neighborhood of a polynomial $f$ is denoted $U_{t,\epsilon}(f)$, where $1/t>M(f)$.  

\begin{lemma} \label{B intersection}
For any $f\in \MP_d$ and $t>0$ such that $1/t > d\cdot M(f)$, we have 
  $$\cB(f,t) = \bigcap_{\eps>0} \pi^{-1} U_{t,\eps}(f).$$
\end{lemma}

\proof
The set $\cB(f,t)$ is clearly contained in the nested intersection, because a conformal conjugacy to $f|X_t(f)$ is an isometry with respect to the conformal metric $|\omega|$.  We now prove the other inclusion.  Any polynomial $g$ in $\cap_{\epsilon>0}\pi^{-1} U_{t,\eps}(f)$  is, on $\{t < G_g < \frac{1}{d\cdot t}\}$, isometrically conjugate to $f$ on $\{t < G_f < \frac{1}{d\cdot t}\}$.  The condition on $t$ guarantees that higher up on the domains $X_{\frac{1}{d\cdot t}}(f)$, $X_{\frac{1}{d\cdot t}}(g)$, the maps $f$ and $g$ are ramified only at the point at infinity.  It follows that this conjugacy extends uniquely to an isometric, hence holomorphic, conjugacy $X_t(f) \to X_t(g)$. So $g \in \cB(f,t)$.
\qed

\subsection{Completing the proof that $\pi$ has connected fibers}
Below, we say that a value $t>0$ is {\em generic} for $f$ if the grand orbits of the critical points do not intersect $\{G_f=t\}$.

\begin{lemma} \label{B connected}
For every $f$ and each generic value $t$ such that $0 < t < \frac{1}{d\cdot M(f)}$, the set $\cB(f,t)$ is connected.
\end{lemma}

\proof
Fix $f$ and a generic value of $t$ with $0< t < \frac{1}{d\cdot M(f)}$.  

Fix $f_1\in \cB(f,t)$, and let $U_\eps \subset \cM_d$ be the connected component of $\pi^{-1}U_{t,\eps}(f)$ containing $f_1$.  We will show that $f_1$ can be connected by a path in $U_\eps$ to $\cS(f,t)$.  Because $\cS(f,t)$ is connected (Corollary \ref{S is path-connected}), it follows that $\cB(f,t)$ is contained in the connected set $U_\eps$.  From Lemma \ref{B intersection}, we have 
	$$\cB(f,t) = \bigcap_{\eps>0} U_\eps$$
and therefore $\cB(f,t)$ is connected.

We first recall the concept of an active critical point; see \cite[\S 4.1]{McMullen:CDR}, \cite[\S 2.1]{Dujardin:Favre:critical}.  Equation (\ref{eqn:param}) gives a map $\rho: \mathcal{H} \times \C \to \cP_d$ parameterizing polynomials by the locations of critical points and constant term.   The $j$th  critical point of a polynomial is {\em active} at a parameter $({\bf c}, a_0) \in \mathcal{H}\times\C$ if the sequence of analytic maps $(c_1, \ldots, c_{d-1}; a) \mapsto \rho(c_1, \ldots, c_{d-1}; a)^{\circ n}(c_j)$, $n \in \N$, fails to be a normal family at $({\bf c}, a_0)$.  If a polynomial has an active critical point, then there exist arbitrarily small perturbations for which this critical point escapes to $\infty$ under iteration.  The locus of polynomials with an active critical point is the {\em bifurcation locus}.  If a polynomial with $m$ escaping critical points has an active critical point, then there exist arbitrarily small perturbations with strictly greater than $m$ escaping critical points.   

In this paragraph, we prove that there is a path $f_s, s \in [0,1]$ contained in $U_\eps$  joining $f_1$ to a map $f_0$ in the shift locus.  If $f_1$ lies in the closure of the shift locus, this is clear.  Otherwise, by the density of structurally stable maps in the family $\MP_d$ (\cite[Cor. 2.8]{McS:QCIII}), there exist arbitrarily small perturbations of $f_1$ which are structurally stable.  
Using  quasiconformal deformations supported on the filled-in Julia set of $f_1$, one finds a path of polynomials converging to a polynomial $f_2$ in the bifurcation locus; the arguments are identical to those given in \cite[\S 5]{DP:hausdorffization}.  An arbitrarily small perturbation of this latter polynomial increases the number of escaping critical points.  By induction, we construct the desired path from $f_1$ to a polynomial $f_0$ in the shift locus.  
The path $f_s$ so constructed is obtained via a sequence of two kinds of modifications: (i) 
arbitrarily small perturbations, and (ii) quasiconformal deformations which do not affect the basin of infinity.  By the continuity of the projection $\pi$, we may assume that this path lies in $U_\epsilon$.  

We now argue that we may assume the path above has in addition the stability property in the hypothesis of Theorem \ref{thm:proved_by_gluing}(3).  On the locus of pairs $(F,z) \in \MP_d \times \C$ for which $\{G_F(z)>0\}$, the map $(z, F) \mapsto G_f(z)$ is pluriharmonic, hence smooth.  The genericity assumption on the given height $t$ thus implies that for the given polynomial $f$, the level sets $\{G_f=t\}$ are nonsingular, and that they remain nonsingular as $s$ varies along a suitably small path.  

Applying Lemma \ref{paths}, we obtain a path in $B(f_0,t)$ joining $f_0$ to an element $g_0$ of $\cS(f_0,t)$.  Since the dynamics above height $t$ along this path is constant, this path lies in $U_\epsilon$.  

Applying Theorem \ref{thm:proved_by_gluing}(3), we obtain a path $s \mapsto g_s \in \cS_d\cap\{g: m(g) \geq t\}$ satisfying $\pi_t(g_s) = \pi_t(f_s)$ for all $s \in [0,1]$.  Since the dynamics above height $t$ along this path is constant, this path too lies in $U_\epsilon$.  By construction, $g_1 \in \cS(f,t)$.  
\qed

\bigskip\noindent
{\em Proof of Theorem \ref{pi}.}  Continuity and properness of $\pi: \MP_d \to \cB_d$ are included in the statement of Lemma \ref{continuous}.  For each point $(X(f), f)$ in $\cB_d$, its fiber is exactly 
	$$\pi^{-1}((X(f), f)) = \bigcap_{\mbox{\tiny generic }\, t>0} \cB(f,t) = \bigcap_{t>0} \cB(f,t),$$
because the sets $\cB(f,t)$ are nested and generic $t$ are dense.
For generic $t$ small enough, the set $\cB(f,t)$ is connected by Lemma \ref{B connected}; therefore $\pi^{-1}((X(f),f))$ is connected.  Finally, Lemma \ref{homeo} states that $\pi$ is a homeomorphism on the shift locus.
\qed

\medskip
Though the fibers of $\pi$ are connected, our methods do not show that they are path-connected.  For example, it is not known if the Mandelbrot set is path-connected. 

\bigskip\bigskip

\def\cprime{$'$}

\vskip 0.5in

\noindent
\textsc{\small Laura G. DeMarco, Department of Mathematics,  University of Illinois at Chicago,    demarco@math.uic.edu }\\

\noindent
\textsc{\small Kevin M. Pilgrim, Department of Mathematics, Indiana University, \\ pilgrim@indiana.edu}

 \end{document}